\documentstyle[11pt]{article}
\newtheorem{theorem}{Theorem}[section]
\newtheorem{lemma}[theorem]{Lemma}
\newtheorem{corollary}[theorem]{Corollary}

\newtheorem{proposition}[theorem]{Proposition}

\makeindex \makeglossary
\begin{document}
%
%

\long\def\ig#1{\relax}
\ig{Thanks to Roberto Minio for this def'n.  Compare the def'n of
\comment in AMSTeX.}

\newcount \coefa
\newcount \coefb
\newcount \coefc
\newcount\tempcounta
\newcount\tempcountb
\newcount\tempcountc
\newcount\tempcountd
\newcount\xext
\newcount\yext
\newcount\xoff
\newcount\yoff
\newcount\gap%
\newcount\arrowtypea
\newcount\arrowtypeb
\newcount\arrowtypec
\newcount\arrowtyped
\newcount\arrowtypee
\newcount\height
\newcount\width
\newcount\xpos
\newcount\ypos
\newcount\run
\newcount\rise
\newcount\arrowlength
\newcount\halflength
\newcount\arrowtype
\newdimen\tempdimen
\newdimen\xlen
\newdimen\ylen
\newsavebox{\tempboxa}%
\newsavebox{\tempboxb}%
\newsavebox{\tempboxc}%

\makeatletter
\setlength{\unitlength}{.01em}%
\def\settypes(#1,#2,#3){\arrowtypea#1 \arrowtypeb#2 \arrowtypec#3}
\def\settoheight#1#2{\setbox\@tempboxa\hbox{#2}#1\ht\@tempboxa\relax}%
\def\settodepth#1#2{\setbox\@tempboxa\hbox{#2}#1\dp\@tempboxa\relax}%
\def\settokens[#1`#2`#3`#4]{%
     \def\tokena{#1}\def\tokenb{#2}\def\tokenc{#3}\def\tokend{#4}}
\def\setsqparms[#1`#2`#3`#4;#5`#6]{%
\arrowtypea #1
\arrowtypeb #2
\arrowtypec #3
\arrowtyped #4
\width #5
\height #6
}
\def\setpos(#1,#2){\xpos=#1 \ypos#2}

\def\bfig{\begin{picture}(\xext,\yext)(\xoff,\yoff)}
\def\efig{\end{picture}}

\def\putbox(#1,#2)#3{\put(#1,#2){\makebox(0,0){$#3$}}}

\def\settriparms[#1`#2`#3;#4]{\settripairparms[#1`#2`#3`1`1;#4]}%

\def\settripairparms[#1`#2`#3`#4`#5;#6]{%
\arrowtypea #1
\arrowtypeb #2
\arrowtypec #3
\arrowtyped #4
\arrowtypee #5
\width #6
\height #6
}

\def\resetparms{\settripairparms[1`1`1`1`1;500]\width 500}

\resetparms

\def\mvector(#1,#2)#3{
\put(0,0){\vector(#1,#2){#3}}%
\put(0,0){\vector(#1,#2){30}}%
}
\def\evector(#1,#2)#3{{
\arrowlength #3
\put(0,0){\vector(#1,#2){\arrowlength}}%
\advance \arrowlength by-30
\put(0,0){\vector(#1,#2){\arrowlength}}%
}}

\def\horsize#1#2{%
\settowidth{\tempdimen}{$#2$}%
#1=\tempdimen
\divide #1 by\unitlength
}

\def\vertsize#1#2{%
\settoheight{\tempdimen}{$#2$}%
#1=\tempdimen
\settodepth{\tempdimen}{$#2$}%
\advance #1 by\tempdimen
\divide #1 by\unitlength
}

\def\vertadjust[#1`#2`#3]{%
\vertsize{\tempcounta}{#1}%
\vertsize{\tempcountb}{#2}%
\ifnum \tempcounta<\tempcountb \tempcounta=\tempcountb \fi
\divide\tempcounta by2
\vertsize{\tempcountb}{#3}%
\ifnum \tempcountb>0 \advance \tempcountb by20 \fi
\ifnum \tempcounta<\tempcountb \tempcounta=\tempcountb \fi
}

\def\horadjust[#1`#2`#3]{%
\horsize{\tempcounta}{#1}%
\horsize{\tempcountb}{#2}%
\ifnum \tempcounta<\tempcountb \tempcounta=\tempcountb \fi
\divide\tempcounta by20
\horsize{\tempcountb}{#3}%
\ifnum \tempcountb>0 \advance \tempcountb by60 \fi
\ifnum \tempcounta<\tempcountb \tempcounta=\tempcountb \fi
}

\ig{ In this procedure, #1 is the paramater that sticks out all the way,
#2 sticks out the least and #3 is a label sticking out half way.  #4 is
the amount of the offset.}

\def\sladjust[#1`#2`#3]#4{%
\tempcountc=#4
\horsize{\tempcounta}{#1}%
\divide \tempcounta by2
\horsize{\tempcountb}{#2}%
\divide \tempcountb by2
\advance \tempcountb by-\tempcountc
\ifnum \tempcounta<\tempcountb \tempcounta=\tempcountb\fi
\divide \tempcountc by2
\horsize{\tempcountb}{#3}%
\advance \tempcountb by-\tempcountc
\ifnum \tempcountb>0 \advance \tempcountb by80\fi
\ifnum \tempcounta<\tempcountb \tempcounta=\tempcountb\fi
\advance\tempcounta by20
}

\def\putvector(#1,#2)(#3,#4)#5#6{{%
\xpos=#1
\ypos=#2
\run=#3
\rise=#4
\arrowlength=#5
\arrowtype=#6
\ifnum \arrowtype<0
    \ifnum \run=0
        \advance \ypos by-\arrowlength
    \else
        \tempcounta \arrowlength
        \multiply \tempcounta by\rise
        \divide \tempcounta by\run
        \ifnum\run>0
            \advance \xpos by\arrowlength
            \advance \ypos by\tempcounta
        \else
            \advance \xpos by-\arrowlength
            \advance \ypos by-\tempcounta
        \fi
    \fi
    \multiply \arrowtype by-1
    \multiply \rise by-1
    \multiply \run by-1
\fi
\ifnum \arrowtype=1
    \put(\xpos,\ypos){\vector(\run,\rise){\arrowlength}}%
\else\ifnum \arrowtype=2
    \put(\xpos,\ypos){\mvector(\run,\rise)\arrowlength}%
\else\ifnum\arrowtype=3
    \put(\xpos,\ypos){\evector(\run,\rise){\arrowlength}}%
\fi\fi\fi
}}

\def\putsplitvector(#1,#2)#3#4{
\xpos #1
\ypos #2
\arrowtype #4
\halflength #3
\arrowlength #3
\gap 140
\advance \halflength by-\gap
\divide \halflength by2
\ifnum \arrowtype=1
    \put(\xpos,\ypos){\line(0,-1){\halflength}}%
    \advance\ypos by-\halflength
    \advance\ypos by-\gap
    \put(\xpos,\ypos){\vector(0,-1){\halflength}}%
\else\ifnum \arrowtype=2
    \put(\xpos,\ypos){\line(0,-1)\halflength}%
    \put(\xpos,\ypos){\vector(0,-1)3}%
    \advance\ypos by-\halflength
    \advance\ypos by-\gap
    \put(\xpos,\ypos){\vector(0,-1){\halflength}}%
\else\ifnum\arrowtype=3
    \put(\xpos,\ypos){\line(0,-1)\halflength}%
    \advance\ypos by-\halflength
    \advance\ypos by-\gap
    \put(\xpos,\ypos){\evector(0,-1){\halflength}}%
\else\ifnum \arrowtype=-1
    \advance \ypos by-\arrowlength
    \put(\xpos,\ypos){\line(0,1){\halflength}}%
    \advance\ypos by\halflength
    \advance\ypos by\gap
    \put(\xpos,\ypos){\vector(0,1){\halflength}}%
\else\ifnum \arrowtype=-2
    \advance \ypos by-\arrowlength
    \put(\xpos,\ypos){\line(0,1)\halflength}%
    \put(\xpos,\ypos){\vector(0,1)3}%
    \advance\ypos by\halflength
    \advance\ypos by\gap
    \put(\xpos,\ypos){\vector(0,1){\halflength}}%
\else\ifnum\arrowtype=-3
    \advance \ypos by-\arrowlength
    \put(\xpos,\ypos){\line(0,1)\halflength}%
    \advance\ypos by\halflength
    \advance\ypos by\gap
    \put(\xpos,\ypos){\evector(0,1){\halflength}}%
\fi\fi\fi\fi\fi\fi
}

\def\putmorphism(#1)(#2,#3)[#4`#5`#6]#7#8#9{{%
\run #2
\rise #3
\ifnum\rise=0
  \puthmorphism(#1)[#4`#5`#6]{#7}{#8}{#9}%
\else\ifnum\run=0
  \putvmorphism(#1)[#4`#5`#6]{#7}{#8}{#9}%
\else
\setpos(#1)%
\arrowlength #7
\arrowtype #8
\ifnum\run=0
\else\ifnum\rise=0
\else
\ifnum\run>0
    \coefa=1
\else
   \coefa=-1
\fi
\ifnum\arrowtype>0
   \coefb=0
   \coefc=-1
\else
   \coefb=\coefa
   \coefc=1
   \arrowtype=-\arrowtype
\fi
\width=2
\multiply \width by\run
\divide \width by\rise
\ifnum \width<0  \width=-\width\fi
\advance\width by60
\if l#9 \width=-\width\fi
\putbox(\xpos,\ypos){#4}
{\multiply \coefa by\arrowlength
\advance\xpos by\coefa
\multiply \coefa by\rise
\divide \coefa by\run
\advance \ypos by\coefa
\putbox(\xpos,\ypos){#5} }%
{\multiply \coefa by\arrowlength
\divide \coefa by2
\advance \xpos by\coefa
\advance \xpos by\width
\multiply \coefa by\rise
\divide \coefa by\run
\advance \ypos by\coefa
\if l#9%
   \put(\xpos,\ypos){\makebox(0,0)[r]{$#6$}}%
\else\if r#9%
   \put(\xpos,\ypos){\makebox(0,0)[l]{$#6$}}%
\fi\fi }%
{\multiply \rise by-\coefc
\multiply \run by-\coefc
\multiply \coefb by\arrowlength
\advance \xpos by\coefb
\multiply \coefb by\rise
\divide \coefb by\run
\advance \ypos by\coefb
\multiply \coefc by70
\advance \ypos by\coefc
\multiply \coefc by\run
\divide \coefc by\rise
\advance \xpos by\coefc
\multiply \coefa by140
\multiply \coefa by\run
\divide \coefa by\rise
\advance \arrowlength by\coefa
\ifnum \arrowtype=1
   \put(\xpos,\ypos){\vector(\run,\rise){\arrowlength}}%
\else\ifnum\arrowtype=2
   \put(\xpos,\ypos){\mvector(\run,\rise){\arrowlength}}%
\else\ifnum\arrowtype=3
   \put(\xpos,\ypos){\evector(\run,\rise){\arrowlength}}%
\fi\fi\fi}\fi\fi\fi\fi}}

\def\puthmorphism(#1,#2)[#3`#4`#5]#6#7#8{{%
\xpos #1
\ypos #2
\width #6
\arrowlength #6
\putbox(\xpos,\ypos){#3\vphantom{#4}}%
{\advance \xpos by\arrowlength
\putbox(\xpos,\ypos){\vphantom{#3}#4}}%
\horsize{\tempcounta}{#3}%
\horsize{\tempcountb}{#4}%
\divide \tempcounta by2
\divide \tempcountb by2
\advance \tempcounta by30
\advance \tempcountb by30
\advance \xpos by\tempcounta
\advance \arrowlength by-\tempcounta
\advance \arrowlength by-\tempcountb
\putvector(\xpos,\ypos)(1,0){\arrowlength}{#7}%
\divide \arrowlength by2
\advance \xpos by\arrowlength
\vertsize{\tempcounta}{#5}%
\divide\tempcounta by2
\advance \tempcounta by20
\if a#8 %
   \advance \ypos by\tempcounta
   \putbox(\xpos,\ypos){#5}%
\else
   \advance \ypos by-\tempcounta
   \putbox(\xpos,\ypos){#5}%
\fi}}

\def\putvmorphism(#1,#2)[#3`#4`#5]#6#7#8{{%
\xpos #1
\ypos #2
\arrowlength #6
\arrowtype #7
\settowidth{\xlen}{$#5$}%
\putbox(\xpos,\ypos){#3}%
{\advance \ypos by-\arrowlength
\putbox(\xpos,\ypos){#4}}%
{\advance\arrowlength by-140
\advance \ypos by-70
\ifdim\xlen>0pt
   \if m#8%
      \putsplitvector(\xpos,\ypos){\arrowlength}{\arrowtype}%
   \else
      \putvector(\xpos,\ypos)(0,-1){\arrowlength}{\arrowtype}%
   \fi
\else
   \putvector(\xpos,\ypos)(0,-1){\arrowlength}{\arrowtype}%
\fi}%
\ifdim\xlen>0pt
   \divide \arrowlength by2
   \advance\ypos by-\arrowlength
   \if l#8%
      \advance \xpos by-40
      \put(\xpos,\ypos){\makebox(0,0)[r]{$#5$}}%
   \else\if r#8%
      \advance \xpos by40
      \put(\xpos,\ypos){\makebox(0,0)[l]{$#5$}}%
   \else
      \putbox(\xpos,\ypos){#5}%
   \fi\fi
\fi
}}

\def\topadjust[#1`#2`#3]{%
\yoff=10
\vertadjust[#1`#2`{#3}]%
\advance \yext by\tempcounta
\advance \yext by 10
}
\def\botadjust[#1`#2`#3]{%
\vertadjust[#1`#2`{#3}]%
\advance \yext by\tempcounta
\advance \yoff by-\tempcounta
}
\def\leftadjust[#1`#2`#3]{%
\xoff=0
\horadjust[#1`#2`{#3}]%
\advance \xext by\tempcounta
\advance \xoff by-\tempcounta
}
\def\rightadjust[#1`#2`#3]{%
\horadjust[#1`#2`{#3}]%
\advance \xext by\tempcounta
}
\def\rightsladjust[#1`#2`#3]{%
\sladjust[#1`#2`{#3}]{\width}%
\advance \xext by\tempcounta
}
\def\leftsladjust[#1`#2`#3]{%
\xoff=0
\sladjust[#1`#2`{#3}]{\width}%
\advance \xext by\tempcounta
\advance \xoff by-\tempcounta
}
\def\adjust[#1`#2;#3`#4;#5`#6;#7`#8]{%
\topadjust[#1``{#2}]
\leftadjust[#3``{#4}]
\rightadjust[#5``{#6}]
\botadjust[#7``{#8}]}

\def\putsquarep<#1>(#2)[#3;#4`#5`#6`#7]{{%
\setsqparms[#1]%
\setpos(#2)%
\settokens[#3]%
\puthmorphism(\xpos,\ypos)[\tokenc`\tokend`{#7}]{\width}{\arrowtyped}b%
\advance\ypos by \height
\puthmorphism(\xpos,\ypos)[\tokena`\tokenb`{#4}]{\width}{\arrowtypea}a%
\putvmorphism(\xpos,\ypos)[``{#5}]{\height}{\arrowtypeb}l%
\advance\xpos by \width
\putvmorphism(\xpos,\ypos)[``{#6}]{\height}{\arrowtypec}r%
}}

\def\putsquare{\@ifnextchar <{\putsquarep}{\putsquarep%
   <\arrowtypea`\arrowtypeb`\arrowtypec`\arrowtyped;\width`\height>}}
\def\square{\@ifnextchar< {\squarep}{\squarep
   <\arrowtypea`\arrowtypeb`\arrowtypec`\arrowtyped;\width`\height>}}
\def\squarep<#1>[#2`#3`#4`#5;#6`#7`#8`#9]{{
\setsqparms[#1]
\xext=\width                                          
\yext=\height                                         
\topadjust[#2`#3`{#6}]
\botadjust[#4`#5`{#9}]
\leftadjust[#2`#4`{#7}]
\rightadjust[#3`#5`{#8}]
\begin{picture}(\xext,\yext)(\xoff,\yoff)
\putsquarep<\arrowtypea`\arrowtypeb`\arrowtypec`\arrowtyped;\width`\height>%
(0,0)[#2`#3`#4`#5;#6`#7`#8`{#9}]%
\end{picture}%
}}

\def\putptrianglep<#1>(#2,#3)[#4`#5`#6;#7`#8`#9]{{%
\settriparms[#1]%
\xpos=#2 \ypos=#3
\advance\ypos by \height
\puthmorphism(\xpos,\ypos)[#4`#5`{#7}]{\height}{\arrowtypea}a%
\putvmorphism(\xpos,\ypos)[`#6`{#8}]{\height}{\arrowtypeb}l%
\advance\xpos by\height
\putmorphism(\xpos,\ypos)(-1,-1)[``{#9}]{\height}{\arrowtypec}r%
}}

\def\putptriangle{\@ifnextchar <{\putptrianglep}{\putptrianglep
   <\arrowtypea`\arrowtypeb`\arrowtypec;\height>}}
\def\ptriangle{\@ifnextchar <{\ptrianglep}{\ptrianglep
   <\arrowtypea`\arrowtypeb`\arrowtypec;\height>}}

\def\ptrianglep<#1>[#2`#3`#4;#5`#6`#7]{{
\settriparms[#1]%
\width=\height                         
\xext=\width                           
\yext=\width                           
\topadjust[#2`#3`{#5}]
\botadjust[#3``]
\leftadjust[#2`#4`{#6}]
\rightsladjust[#3`#4`{#7}]
\begin{picture}(\xext,\yext)(\xoff,\yoff)
\putptrianglep<\arrowtypea`\arrowtypeb`\arrowtypec;\height>%
(0,0)[#2`#3`#4;#5`#6`{#7}]%
\end{picture}%
}}

\def\putqtrianglep<#1>(#2,#3)[#4`#5`#6;#7`#8`#9]{{%
\settriparms[#1]%
\xpos=#2 \ypos=#3
\advance\ypos by\height
\puthmorphism(\xpos,\ypos)[#4`#5`{#7}]{\height}{\arrowtypea}a%
\putmorphism(\xpos,\ypos)(1,-1)[``{#8}]{\height}{\arrowtypeb}l%
\advance\xpos by\height
\putvmorphism(\xpos,\ypos)[`#6`{#9}]{\height}{\arrowtypec}r%
}}

\def\putqtriangle{\@ifnextchar <{\putqtrianglep}{\putqtrianglep
   <\arrowtypea`\arrowtypeb`\arrowtypec;\height>}}
\def\qtriangle{\@ifnextchar <{\qtrianglep}{\qtrianglep
   <\arrowtypea`\arrowtypeb`\arrowtypec;\height>}}

\def\qtrianglep<#1>[#2`#3`#4;#5`#6`#7]{{
\settriparms[#1]
\width=\height                         
\xext=\width                           
\yext=\height                          
\topadjust[#2`#3`{#5}]
\botadjust[#4``]
\leftsladjust[#2`#4`{#6}]
\rightadjust[#3`#4`{#7}]
\begin{picture}(\xext,\yext)(\xoff,\yoff)
\putqtrianglep<\arrowtypea`\arrowtypeb`\arrowtypec;\height>%
(0,0)[#2`#3`#4;#5`#6`{#7}]%
\end{picture}%
}}

\def\putdtrianglep<#1>(#2,#3)[#4`#5`#6;#7`#8`#9]{{%
\settriparms[#1]%
\xpos=#2 \ypos=#3
\puthmorphism(\xpos,\ypos)[#5`#6`{#9}]{\height}{\arrowtypec}b%
\advance\xpos by \height \advance\ypos by\height
\putmorphism(\xpos,\ypos)(-1,-1)[``{#7}]{\height}{\arrowtypea}l%
\putvmorphism(\xpos,\ypos)[#4``{#8}]{\height}{\arrowtypeb}r%
}}

\def\putdtriangle{\@ifnextchar <{\putdtrianglep}{\putdtrianglep
   <\arrowtypea`\arrowtypeb`\arrowtypec;\height>}}
\def\dtriangle{\@ifnextchar <{\dtrianglep}{\dtrianglep
   <\arrowtypea`\arrowtypeb`\arrowtypec;\height>}}

\def\dtrianglep<#1>[#2`#3`#4;#5`#6`#7]{{
\settriparms[#1]
\width=\height                         
\xext=\width                           
\yext=\height                          
\topadjust[#2``]
\botadjust[#3`#4`{#7}]
\leftsladjust[#3`#2`{#5}]
\rightadjust[#2`#4`{#6}]
\begin{picture}(\xext,\yext)(\xoff,\yoff)
\putdtrianglep<\arrowtypea`\arrowtypeb`\arrowtypec;\height>%
(0,0)[#2`#3`#4;#5`#6`{#7}]%
\end{picture}%
}}

\def\putbtrianglep<#1>(#2,#3)[#4`#5`#6;#7`#8`#9]{{%
\settriparms[#1]%
\xpos=#2 \ypos=#3
\puthmorphism(\xpos,\ypos)[#5`#6`{#9}]{\height}{\arrowtypec}b%
\advance\ypos by\height
\putmorphism(\xpos,\ypos)(1,-1)[``{#8}]{\height}{\arrowtypeb}r%
\putvmorphism(\xpos,\ypos)[#4``{#7}]{\height}{\arrowtypea}l%
}}

\def\putbtriangle{\@ifnextchar <{\putbtrianglep}{\putbtrianglep
   <\arrowtypea`\arrowtypeb`\arrowtypec;\height>}}
\def\btriangle{\@ifnextchar <{\btrianglep}{\btrianglep
   <\arrowtypea`\arrowtypeb`\arrowtypec;\height>}}

\def\btrianglep<#1>[#2`#3`#4;#5`#6`#7]{{
\settriparms[#1]
\width=\height                         
\xext=\width                           
\yext=\height                          
\topadjust[#2``]
\botadjust[#3`#4`{#7}]
\leftadjust[#2`#3`{#5}]
\rightsladjust[#4`#2`{#6}]
\begin{picture}(\xext,\yext)(\xoff,\yoff)
\putbtrianglep<\arrowtypea`\arrowtypeb`\arrowtypec;\height>%
(0,0)[#2`#3`#4;#5`#6`{#7}]%
\end{picture}%
}}

\def\putAtrianglep<#1>(#2,#3)[#4`#5`#6;#7`#8`#9]{{%
\settriparms[#1]%
\xpos=#2 \ypos=#3
{\multiply \height by2
\puthmorphism(\xpos,\ypos)[#5`#6`{#9}]{\height}{\arrowtypec}b}%
\advance\xpos by\height \advance\ypos by\height
\putmorphism(\xpos,\ypos)(-1,-1)[#4``{#7}]{\height}{\arrowtypea}l%
\putmorphism(\xpos,\ypos)(1,-1)[``{#8}]{\height}{\arrowtypeb}r%
}}

\def\putAtriangle{\@ifnextchar <{\putAtrianglep}{\putAtrianglep
   <\arrowtypea`\arrowtypeb`\arrowtypec;\height>}}
\def\Atriangle{\@ifnextchar <{\Atrianglep}{\Atrianglep
   <\arrowtypea`\arrowtypeb`\arrowtypec;\height>}}

\def\Atrianglep<#1>[#2`#3`#4;#5`#6`#7]{{
\settriparms[#1]
\width=\height                         
\xext=\width                           
\yext=\height                          
\topadjust[#2``]
\botadjust[#3`#4`{#7}]
\multiply \xext by2 
\leftsladjust[#3`#2`{#5}]
\rightsladjust[#4`#2`{#6}]
\begin{picture}(\xext,\yext)(\xoff,\yoff)%
\putAtrianglep<\arrowtypea`\arrowtypeb`\arrowtypec;\height>%
(0,0)[#2`#3`#4;#5`#6`{#7}]%
\end{picture}%
}}

\def\putAtrianglepairp<#1>(#2)[#3;#4`#5`#6`#7`#8]{{
\settripairparms[#1]%
\setpos(#2)%
\settokens[#3]%
\puthmorphism(\xpos,\ypos)[\tokenb`\tokenc`{#7}]{\height}{\arrowtyped}b%
\advance\xpos by\height
\advance\ypos by\height
\putmorphism(\xpos,\ypos)(-1,-1)[\tokena``{#4}]{\height}{\arrowtypea}l%
\putvmorphism(\xpos,\ypos)[``{#5}]{\height}{\arrowtypeb}m%
\putmorphism(\xpos,\ypos)(1,-1)[``{#6}]{\height}{\arrowtypec}r%
}}

\def\putAtrianglepair{\@ifnextchar <{\putAtrianglepairp}{\putAtrianglepairp%
   <\arrowtypea`\arrowtypeb`\arrowtypec`\arrowtyped`\arrowtypee;\height>}}
\def\Atrianglepair{\@ifnextchar <{\Atrianglepairp}{\Atrianglepairp%
   <\arrowtypea`\arrowtypeb`\arrowtypec`\arrowtyped`\arrowtypee;\height>}}

\def\Atrianglepairp<#1>[#2;#3`#4`#5`#6`#7]{{%
\settripairparms[#1]%
\settokens[#2]%
\width=\height
\xext=\width
\yext=\height
\topadjust[\tokena``]%
\vertadjust[\tokenb`\tokenc`{#6}]
\tempcountd=\tempcounta                       
\vertadjust[\tokenc`\tokend`{#7}]
\ifnum\tempcounta<\tempcountd                 
\tempcounta=\tempcountd\fi                    
\advance \yext by\tempcounta                  
\advance \yoff by-\tempcounta                 %
\multiply \xext by2 
\leftsladjust[\tokenb`\tokena`{#3}]
\rightsladjust[\tokend`\tokena`{#5}]%
\begin{picture}(\xext,\yext)(\xoff,\yoff)%
\putAtrianglepairp
<\arrowtypea`\arrowtypeb`\arrowtypec`\arrowtyped`\arrowtypee;\height>%
(0,0)[#2;#3`#4`#5`#6`{#7}]%
\end{picture}%
}}

\def\putVtrianglep<#1>(#2,#3)[#4`#5`#6;#7`#8`#9]{{%
\settriparms[#1]%
\xpos=#2 \ypos=#3
\advance\ypos by\height
{\multiply\height by2
\puthmorphism(\xpos,\ypos)[#4`#5`{#7}]{\height}{\arrowtypea}a}%
\putmorphism(\xpos,\ypos)(1,-1)[`#6`{#8}]{\height}{\arrowtypeb}l%
\advance\xpos by\height
\advance\xpos by\height
\putmorphism(\xpos,\ypos)(-1,-1)[``{#9}]{\height}{\arrowtypec}r%
}}

\def\putVtriangle{\@ifnextchar <{\putVtrianglep}{\putVtrianglep
   <\arrowtypea`\arrowtypeb`\arrowtypec;\height>}}
\def\Vtriangle{\@ifnextchar <{\Vtrianglep}{\Vtrianglep
   <\arrowtypea`\arrowtypeb`\arrowtypec;\height>}}

\def\Vtrianglep<#1>[#2`#3`#4;#5`#6`#7]{{
\settriparms[#1]
\width=\height                         
\xext=\width                           
\yext=\height                          
\topadjust[#2`#3`{#5}]
\botadjust[#4``]
\multiply \xext by2 
\leftsladjust[#2`#3`{#6}]
\rightsladjust[#3`#4`{#7}]
\begin{picture}(\xext,\yext)(\xoff,\yoff)%
\putVtrianglep<\arrowtypea`\arrowtypeb`\arrowtypec;\height>%
(0,0)[#2`#3`#4;#5`#6`{#7}]%
\end{picture}%
}}

\def\putVtrianglepairp<#1>(#2)[#3;#4`#5`#6`#7`#8]{{
\settripairparms[#1]%
\setpos(#2)%
\settokens[#3]%
\advance\ypos by\height
\putmorphism(\xpos,\ypos)(1,-1)[`\tokend`{#6}]{\height}{\arrowtypec}l%
\puthmorphism(\xpos,\ypos)[\tokena`\tokenb`{#4}]{\height}{\arrowtypea}a%
\advance\xpos by\height
\putvmorphism(\xpos,\ypos)[``{#7}]{\height}{\arrowtyped}m%
\advance\xpos by\height
\putmorphism(\xpos,\ypos)(-1,-1)[``{#8}]{\height}{\arrowtypee}r%
}}

\def\putVtrianglepair{\@ifnextchar <{\putVtrianglepairp}{\putVtrianglepairp%
    <\arrowtypea`\arrowtypeb`\arrowtypec`\arrowtyped`\arrowtypee;\height>}}
\def\Vtrianglepair{\@ifnextchar <{\Vtrianglepairp}{\Vtrianglepairp%
    <\arrowtypea`\arrowtypeb`\arrowtypec`\arrowtyped`\arrowtypee;\height>}}

\def\Vtrianglepairp<#1>[#2;#3`#4`#5`#6`#7]{{%
\settripairparms[#1]%
\settokens[#2]
\xext=\height                  
\width=\height                 
\yext=\height                  
\vertadjust[\tokena`\tokenb`{#4}]
\tempcountd=\tempcounta        
\vertadjust[\tokenb`\tokenc`{#5}]
\ifnum\tempcounta<\tempcountd%
\tempcounta=\tempcountd\fi
\advance \yext by\tempcounta
\botadjust[\tokend``]%
\multiply \xext by2
\leftsladjust[\tokena`\tokend`{#6}]%
\rightsladjust[\tokenc`\tokend`{#7}]%
\begin{picture}(\xext,\yext)(\xoff,\yoff)%
\putVtrianglepairp
<\arrowtypea`\arrowtypeb`\arrowtypec`\arrowtyped`\arrowtypee;\height>%
(0,0)[#2;#3`#4`#5`#6`{#7}]%
\end{picture}%
}}

\def\putCtrianglep<#1>(#2,#3)[#4`#5`#6;#7`#8`#9]{{%
\settriparms[#1]%
\xpos=#2 \ypos=#3
\advance\ypos by\height
\putmorphism(\xpos,\ypos)(1,-1)[``{#9}]{\height}{\arrowtypec}l%
\advance\xpos by\height
\advance\ypos by\height
\putmorphism(\xpos,\ypos)(-1,-1)[#4`#5`{#7}]{\height}{\arrowtypea}l%
{\multiply\height by 2
\putvmorphism(\xpos,\ypos)[`#6`{#8}]{\height}{\arrowtypeb}r}%
}}

\def\putCtriangle{\@ifnextchar <{\putCtrianglep}{\putCtrianglep
    <\arrowtypea`\arrowtypeb`\arrowtypec;\height>}}
\def\Ctriangle{\@ifnextchar <{\Ctrianglep}{\Ctrianglep
    <\arrowtypea`\arrowtypeb`\arrowtypec;\height>}}

\def\Ctrianglep<#1>[#2`#3`#4;#5`#6`#7]{{
\settriparms[#1]
\width=\height                          
\xext=\width                            
\yext=\height                           
\multiply \yext by2 
\topadjust[#2``]
\botadjust[#4``]
\sladjust[#3`#2`{#5}]{\width}
\tempcountd=\tempcounta                 
\sladjust[#3`#4`{#7}]{\width}
\ifnum \tempcounta<\tempcountd          
\tempcounta=\tempcountd\fi              
\advance \xext by\tempcounta            
\advance \xoff by-\tempcounta           %
\rightadjust[#2`#4`{#6}]
\begin{picture}(\xext,\yext)(\xoff,\yoff)%
\putCtrianglep<\arrowtypea`\arrowtypeb`\arrowtypec;\height>%
(0,0)[#2`#3`#4;#5`#6`{#7}]%
\end{picture}%
}}

\def\putDtrianglep<#1>(#2,#3)[#4`#5`#6;#7`#8`#9]{{%
\settriparms[#1]%
\xpos=#2 \ypos=#3
\advance\xpos by\height \advance\ypos by\height
\putmorphism(\xpos,\ypos)(-1,-1)[``{#9}]{\height}{\arrowtypec}r%
\advance\xpos by-\height \advance\ypos by\height
\putmorphism(\xpos,\ypos)(1,-1)[`#5`{#8}]{\height}{\arrowtypeb}r%
{\multiply\height by 2
\putvmorphism(\xpos,\ypos)[#4`#6`{#7}]{\height}{\arrowtypea}l}%
}}

\def\putDtriangle{\@ifnextchar <{\putDtrianglep}{\putDtrianglep
    <\arrowtypea`\arrowtypeb`\arrowtypec;\height>}}
\def\Dtriangle{\@ifnextchar <{\Dtrianglep}{\Dtrianglep
   <\arrowtypea`\arrowtypeb`\arrowtypec;\height>}}

\def\Dtrianglep<#1>[#2`#3`#4;#5`#6`#7]{{
\settriparms[#1]
\width=\height                         
\xext=\height                          
\yext=\height                          
\multiply \yext by2 
\topadjust[#2``]
\botadjust[#4``]
\leftadjust[#2`#4`{#5}]
\sladjust[#3`#2`{#5}]{\height}
\tempcountd=\tempcountd                
\sladjust[#3`#4`{#7}]{\height}
\ifnum \tempcounta<\tempcountd         
\tempcounta=\tempcountd\fi             
\advance \xext by\tempcounta           %
\begin{picture}(\xext,\yext)(\xoff,\yoff)
\putDtrianglep<\arrowtypea`\arrowtypeb`\arrowtypec;\height>%
(0,0)[#2`#3`#4;#5`#6`{#7}]%
\end{picture}%
}}

\def\setrecparms[#1`#2]{\width=#1 \height=#2}%
%

\def\recursep<#1`#2>[#3;#4`#5`#6`#7`#8]{{%
\width=#1 \height=#2
\settokens[#3]
\settowidth{\tempdimen}{$\tokena$}
\ifdim\tempdimen=0pt
  \savebox{\tempboxa}{\hbox{$\tokenb$}}%
  \savebox{\tempboxb}{\hbox{$\tokend$}}%
  \savebox{\tempboxc}{\hbox{$#6$}}%
\else
  \savebox{\tempboxa}{\hbox{$\hbox{$\tokena$}\times\hbox{$\tokenb$}$}}%
  \savebox{\tempboxb}{\hbox{$\hbox{$\tokena$}\times\hbox{$\tokend$}$}}%
  \savebox{\tempboxc}{\hbox{$\hbox{$\tokena$}\times\hbox{$#6$}$}}%
\fi
\ypos=\height
\divide\ypos by 2
\xpos=\ypos
\advance\xpos by \width
\xext=\xpos \yext=\height
\topadjust[#3`\usebox{\tempboxa}`{#4}]%
\botadjust[#5`\usebox{\tempboxb}`{#8}]%
\sladjust[\tokenc`\tokenb`{#5}]{\ypos}%
\tempcountd=\tempcounta
\sladjust[\tokenc`\tokend`{#5}]{\ypos}%
\ifnum \tempcounta<\tempcountd
\tempcounta=\tempcountd\fi
\advance \xext by\tempcounta
\advance \xoff by-\tempcounta
\rightadjust[\usebox{\tempboxa}`\usebox{\tempboxb}`\usebox{\tempboxc}]%
\bfig
\putCtrianglep<-1`1`1;\ypos>(0,0)[`\tokenc`;#5`#6`{#7}]%
\puthmorphism(\ypos,0)[\tokend`\usebox{\tempboxb}`{#8}]{\width}{-1}b%
\puthmorphism(\ypos,\height)[\tokenb`\usebox{\tempboxa}`{#4}]{\width}{-1}a%
\advance\ypos by \width
\putvmorphism(\ypos,\height)[``\usebox{\tempboxc}]{\height}1r%
\efig
}}

\def\recurse{\@ifnextchar <{\recursep}{\recursep<\width`\height>}}

\def\puttwohmorphisms(#1,#2)[#3`#4;#5`#6]#7#8#9{{%
%
\puthmorphism(#1,#2)[#3`#4`]{#7}0a
\ypos=#2
\advance\ypos by 20
\puthmorphism(#1,\ypos)[\phantom{#3}`\phantom{#4}`#5]{#7}{#8}a
\advance\ypos by -40
\puthmorphism(#1,\ypos)[\phantom{#3}`\phantom{#4}`#6]{#7}{#9}b
}}

\def\puttwovmorphisms(#1,#2)[#3`#4;#5`#6]#7#8#9{{%
%
%
%
\putvmorphism(#1,#2)[#3`#4`]{#7}0a
\xpos=#1
\advance\xpos by -20
\putvmorphism(\xpos,#2)[\phantom{#3}`\phantom{#4}`#5]{#7}{#8}l
\advance\xpos by 40
\putvmorphism(\xpos,#2)[\phantom{#3}`\phantom{#4}`#6]{#7}{#9}r
}}

\def\puthcoequalizer(#1)[#2`#3`#4;#5`#6`#7]#8#9{{%
%
\setpos(#1)%
\puttwohmorphisms(\xpos,\ypos)[#2`#3;#5`#6]{#8}11%
\advance\xpos by #8
\puthmorphism(\xpos,\ypos)[\phantom{#3}`#4`#7]{#8}1{#9}
}}

\def\putvcoequalizer(#1)[#2`#3`#4;#5`#6`#7]#8#9{{%
%
%
%
%
\setpos(#1)%
\puttwovmorphisms(\xpos,\ypos)[#2`#3;#5`#6]{#8}11%
\advance\ypos by -#8
\putvmorphism(\xpos,\ypos)[\phantom{#3}`#4`#7]{#8}1{#9}
}}

\def\putthreehmorphisms(#1)[#2`#3;#4`#5`#6]#7(#8)#9{{%
\setpos(#1) \settypes(#8)
\if a#9 %
     \vertsize{\tempcounta}{#5}%
     \vertsize{\tempcountb}{#6}%
     \ifnum \tempcounta<\tempcountb \tempcounta=\tempcountb \fi
\else
     \vertsize{\tempcounta}{#4}%
     \vertsize{\tempcountb}{#5}%
     \ifnum \tempcounta<\tempcountb \tempcounta=\tempcountb \fi
\fi
\advance \tempcounta by 60
\puthmorphism(\xpos,\ypos)[#2`#3`#5]{#7}{\arrowtypeb}{#9}
\advance\ypos by \tempcounta
\puthmorphism(\xpos,\ypos)[\phantom{#2}`\phantom{#3}`#4]{#7}{\arrowtypea}{#9}
\advance\ypos by -\tempcounta \advance\ypos by -\tempcounta
\puthmorphism(\xpos,\ypos)[\phantom{#2}`\phantom{#3}`#6]{#7}{\arrowtypec}{#9}
}}

\def\putarc(#1,#2)[#3`#4`#5]#6#7#8{{%
\xpos #1
\ypos #2
\width #6
\arrowlength #6
\putbox(\xpos,\ypos){#3\vphantom{#4}}%
{\advance \xpos by\arrowlength
\putbox(\xpos,\ypos){\vphantom{#3}#4}}%
\horsize{\tempcounta}{#3}%
\horsize{\tempcountb}{#4}%
\divide \tempcounta by2
\divide \tempcountb by2
\advance \tempcounta by30
\advance \tempcountb by30
\advance \xpos by\tempcounta
\advance \arrowlength by-\tempcounta
\advance \arrowlength by-\tempcountb
\halflength=\arrowlength \divide\halflength by 2
\divide\arrowlength by 5
\put(\xpos,\ypos){\bezier{\arrowlength}(0,0)(50,50)(\halflength,50)}
\ifnum #7=-1 \put(\xpos,\ypos){\vector(-3,-2)0} \fi
\advance\xpos by \halflength
\put(\xpos,\ypos){\xpos=\halflength \advance\xpos by -50
   \bezier{\arrowlength}(0,50)(\xpos,50)(\halflength,0)}
\ifnum #7=1 {\advance \xpos by
   \halflength \put(\xpos,\ypos){\vector(3,-2)0}} \fi
\advance\ypos by 50
\vertsize{\tempcounta}{#5}%
\divide\tempcounta by2
\advance \tempcounta by20
\if a#8 %
   \advance \ypos by\tempcounta
   \putbox(\xpos,\ypos){#5}%
\else
   \advance \ypos by-\tempcounta
   \putbox(\xpos,\ypos){#5}%
\fi
}}

\makeatother

\sloppy

\def\dsum{\sqcup\!\!\!\!\cdot\;}

\newcommand{\lra}{\longrightarrow}
\newcommand{\ra}{\rightarrow}
\newcommand{\lla}{\longleftarrow}
\newcommand{\da}{\downarrow}
\newcommand{\ua}{\uparrow}
\newcommand{\Da}{\Downarrow}
\newcommand{\Ua}{\Uparrow}

\newcommand{\ira}{\longrightarrow\!\!\!\!\!\!\!\!\circ\;\;\,}
\newcommand{\sira}{\rightarrow\!\!\!\!\!\circ\;\,}

\def\phi{\varphi}
\def\o{{\omega}}

\def\bA{{\bf A}}
\def\bM{{\bf M}}
\def\bN{{\bf N}}
\def\bC{{\bf C}}
\def\bP{{\bf P}}
\def\ba{{\bf a}}
\def\bb{{\bf b}}
\def\bc{{\bf c}}
\def\bd{{\bf d}}
\def\bi{{\bf i}}
\def\bm{{\bf m}}
\def\be{{\bf e}}

\def\cB{{\cal B}}
\def\cA{{\cal A}}
\def\cC{{\cal C}}
\def\cD{{\cal D}}
\def\cE{{\cal E}}
\def\cF{{\cal F}}
\def\cG{{\cal G}}
\def\cK{{\cal K}}
\def\cL{{\cal L}}
\def\cN{{\cal N}}
\def\cM{{\cal M}}
\def\cP{{\cal P}}


\def\mlt{{\bf Mlt}}
\def\cat{{\bf Cat}}

\def\oC{{{\omega}Cat}}
\def\oG{{{\omega}Gr}}
\def\mts{{MltSet}}


\def\ofs{{\bf oFs}}
\def\lfs{{\bf lFs}}
\def\ofsl{{\bf oFs}_{loc}}
\def\ofso{{\bf oFs}_{\o}}
\def\ofsm{{\bf oFs}_{\mu}}
\def\pfs{{\bf pFs}}


\def\pComp{{\bf Comp}^{+/1}} 


\def\mComp{{\bf Comp}^{m/1}} 
\def\Comp{{\bf Comp}} 
\def\Shape{{\bf Shape}} 
\def\mShape{{\bf Shape}^{m/1}} 

\def\mnComma{{\bf Comma}^{m/1}_n} 
\def\fs{{\bf Fs}^{+/1}}

\def\mltsets{{\bf MltSets}} 
\def\comp{Comp} 

\pagenumbering{arabic} \setcounter{page}{1}

\title{Multitopes are the same as \\
principal ordered face structures}
\author{Marek Zawadowski
\\
Instytut Matematyki, Uniwersytet Warszawski\\
ul. S.Banacha 2, 00-913 Warszawa, Poland\\
zawado@mimuw.edu.pl\\}

\maketitle
\begin{flushright}
\em Dedicated to Professor F.W.Lawvere\\ on the occasion of his
70th birthday.
\end{flushright}

\begin{abstract}  We show that the category of
principal ordered face structures $\pfs$ is equivalent to the
category of multitopes  $\bf Mlt$. On the way we introduce the notion
of a graded tensor theory to state the abstract
properties of the category of ordered face structures $\ofs$ and show
how $\ofs$ fits into the recent work of T. Leinster and M.
Weber concerning the nerve construction.

MS Classification 18F20, 18C20, 18D10, 18D05 (AMS 2000).
\end{abstract}

\section{Introduction}


In \cite{Z1} the notion of a positive face structure is introduced
and it is shown how it helps to understand the positive-to-one
computads. In \cite{Z2} part of the program of \cite{Z1} was
developed in the many-to-one context, i.e. the notion of an
ordered face structure was introduced and related to the
many-to-one computads. The first part of this paper is a sequel of \cite{Z2}
developing farther part of \cite{Z1} in the many-to-one context.
We show how the category of ordered face structures and monotone maps $\ofs$  can be used to show that the
category of many-to-one computads $\mComp$ is equivalent to the presheaves
category $Set^{\pfs^{op}}$, where $\pfs$ is the full subcategory of $\ofs$ whose
objects are principal ordered face structures. In fact we show
that both categories $\mComp$ and $Set^{\pfs^{op}}$ are equivalent
to the category $Mod_\otimes(\ofs^{op},Set)$ of $Set$-models of the
graded tensor theory $\ofs$. In \cite{HMP} it was shown that
the category $MltSet$ of multitopic set is a presheaf category on the category of multitopes $\bf Mlt$.
In \cite{HMZ} it was shown that the categories $MltSet$ and $\mComp$
are equivalent. Thus as a corollary we get the statement from
the title of this paper.


My main motivation to define the ordered face structures was to have an explicit combinatorial definition of multitopic sets\footnote{Recall that
multitopic category, a weak $\o$-category in the sense of M.Makkai, is a multitopic set with a property, c.f. [M].}
(or what comes to the same the many-to-one computads, c.f. \cite{HMZ})  that allow fairly easy direct manipulation on cells.
For this I wanted to describe not only the shapes of many-to-one indets (=indeterminates) but of all the cells build from them.
To see some pictures and more explanations on this consult introduction to \cite{Z2}.
As there are several other structures that are serving a similar purpose, that I will discuss later, the anonymous referee asked to explain what is the role of the category $\ofs$ and why it is of an interest at all as its definition is not a simple one.
It is not always easy to give a convincing answer to such questions. After my talk  describing the ordered face structures in
Patras (PSSL, April 2008) J. Kock suggested that the recent paper of M. Weber, c.f. \cite{Weber}, could provide a framework for a conceptual explanation what $\ofs$ is. The explanations I will present in the second part of the paper are very much inspired by the work of
T. Leinster \cite{Leinster1} and M. Weber \cite{Weber} but it also goes beyond that. The short answer is that the category $\ofs$ is the category of shapes of all cells, not only indeterminates (=indets), in many-to-one computads. The abstract properties of $\ofs$ are subsumed by the notion of a graded tensor theory.
I can also make a broader but 'non-full' analogy concerning $\ofs$.
It is related to the $\o$-category monad on many-to-one computads in a similar way as the category of simple $\o$-graphs $s\o\cG$ (or globular cardinals) is related to the $\o$-category monad on one-to-one computads $\Comp^{1/1}$, i.e. the free $\o$-categories over $\o$-graphs with morphisms sending indets to indets.
However  the embedding $\ofs \ra \mComp$ is not full and what is even worse it is not full on isomorphisms. The full image of $\ofs$ under this embedding is
the category of ordered face structures and local maps $\ofsl$ which plays the role of the category of many-to-one cardinals in the Leinster-Weber approach. Thus we have here two different categories $\ofs$ and $\ofsl$ where T. Leinster and M. Weber have only one. What I mean by the category of shapes is a bit technical and the precise definition will be given in Section \ref{shape}.


T. Leinster in \cite{Leinster2} explained that he started to love the nerve construction when he discovered
that both category $\Delta$ and the nerve construction (for categories) arise canonically from the free category monad on graphs. This convinced him
that the construction is {\em natural}. Before, he could only acknowledge that the nerve construction is just {\em useful}.
I would consider even two earlier stages in the process of proving 'rights to exists' of a concept.  One, when there is a  {\em construction} of the
object in question which is not totally unrelated to the purpose it serves. Then, if all else failed, there might be  a {\em purity} of style behind the
notion. The reason I explain all this is that I don't have a canonical simple construction that would make T. Leinster believe that the category $\ofs$ can
be naturally derived from \mbox{$\o$-category} monad on the category of many-to-one computads $\mComp$ or possibly some other fundamental construction
related $\mComp$. But I will argue about the three weaker claims.

1. {\em Purity: simple combinatorial data.}  As I mentioned earlier, my main motivation to define the ordered face structures was to have
an explicit combinatorial definition of multitopic sets that allow fairly easy direct manipulation on cells.
I wanted to describe these shapes with the least possible structure.  So in an ordered face structure we have
functions $\gamma$, associating a {\em face} $\gamma(a)$ which is the codomain the cell $a$,  relations $\delta$,
associating a {\em set of faces} $\delta(a)$ that 'constitue' the domain to the cell $a$, and strict order relations
$<^\sim$ that will indicate in case of doubts in what order one should compose the cells.
The structure is kept so simple at the expense of the axioms that are quite involved and do not
look at first sight as something that have much to do with what it was designed for.
To explain how ordered face structures describe many-to-one computads is a long story, see \cite{Z2}.

2. {\em Abstract construction: the category of shapes.}
However there is an abstract definition of the category of shapes that in most considered cases
gives the category which is equivalent to the category of $T$-cardinals considered by T. Leinster and M. Weber
but in the case of many-to-one computads it is equivalent to $\ofs$ rather than the category of cardinals which is in this case $\ofsl$.
This definition of the category of shapes is given in the  section \ref{shape}. It is at least related with the many-to-one computads
from the very beginning but it is rather hard to believe that it might be of any practical use.

3. {\em Usefulness: $\ofs$ generates all the setup of Leinster and Weber.}
In Sections \ref{pra monads} and \ref{setup} I will argue that $\ofs$ is useful as this category alone generates
all the setup it is involved with. That includes the category of many-to-one computads $\mComp$, the $\o$-category monad on the category $\mComp$, the proof that that this monad is a parametric right adjoint, and that $\o$-categories can be considered as some presheaf satisfying an additional condition.


The notion of a graded tensor theory, GT-theory for short, is designed to describe the abstract features of the category $\ofs$.
Any model $M:\cC^{op}\ra \cA$ of an GT-theory $\cC$ in a category $\cA$
gives rise to a functor $\bar{M}:\cA\ra\oC$ from the category $\cA$ to the category of strict $\o$-categories.  This notion was inspired by and
should be compared with the notion of a monoidal globular category, MG-category for short, of M. Batanin, c.f. \cite{Batanin}.
Both notions deal with the $k$-domain and the $k$-codomain operations. Both notions have the $k$-tensor product operations that can
be performed only if the $k$-codomain of the first object agrees with the $k$-domain of the second one. In GT-theories the cylinder operation is not
given explicitly. However there are essential differences.
An GT-theory $\cC$ is a single (rigid) category together with a dimension functor $dim: \cC\ra \bN$ into the linear order of natural numbers $\bN$.
The $k$-tensor operations are required to be functorial, as in MG-category, but the $k$-domain and the $k$-codomain operations are not functorial in general.
Instead all these operations are given together with specified morphisms $\bd^{(k)}_S:\bd^{(k)}_S\ra S$, $\bc^{(k)}_S:\bc^{(k)}_S\ra S$  in $\cC$,
 $\kappa^1_S: S\ra S\otimes_k S'$,  $\kappa^2_{S'}: S\ra S\otimes_k S'$ that explain the relation of the domains  and
the codomains of objects to objects themselves and of the components of the tensor  products to the tensor products. There are isomorphisms relating
these operations as in MG-category but, as an GT-theory is a rigid category, the coherence conditions are satisfied automatically.
Last but not least the category $\ofs$ is a GT-theory but the domain and the codomain operations are not functorial and the truncations of $\ofs$ do not form an MG-category, contrary to a public claim I have made. It is true that the isomorphism classes of objects in an GT-theory can be easily organized into a discrete MG-category.  But this process when applied to $\ofs$ would destroy an essential information about the monotone morphisms. The notion of a model of GT-theory (a functor sending tensor square to pullbacks) is very important in this context but doesn't seem to have an analog in the context of MG-categories.


In Section \ref{pra monads} the setup developed by T. Leinster and M. Weber is recalled but not in the full generality of \cite{Weber} and in
a form that changes the emphasis. So I will not recall the setup here but only point out to the change in the emphasis.
I will discuss only the parametric right adjoint monads on presheaf categories, called pra monads for short. The monadic functor inducing pra monad is called {\em pra monadic}. Among pra functors there are particularly simple ones that arise from factorizations system on small categories. If $(\cE,\cM)$ is a factorization system on a category $\cC$, satisfying a simple additional condition to be found in Section \ref{pra monads}, and $\cC_\cM$ is the category with the objects from $\cC$ and morphisms from $\cM$ then the restriction functor $i^*: Set^{{\cC}^{op}}\lra Set^{{\cC_\cM}^{op}}$ along the inclusion $i:\cC_\cM\ra \cC$ is pra monadic.  Such functors I will call {\em presheaf pra monadic}. In this context the main conclusion of the work of T. Leinster and M. Weber, present in \cite{Weber} in a sightly hidden form, is that any pra monadic functors, arise as pseudo-pullback of a presheaf pra monadic one along a full and faithful functor. Thus it can be thought of as a representation/completeness result for pra monads. In Section \ref{setup} an extension of the above setup is proposed and it is shown how the category $\ofs$ generates all its ingredients.


In the presheaf approaches  to weak categories (as opposed to the algebraic ones) the weak categories are presheaves with some properties.
If we believe that strict $\o$-categories should be special cases of weak ones we need to
find the way how to interpret strict $\o$-categories as appropriate presheaves.
The nerves of $\o$-categories are constructions that do exactly this and (should) provide
abundance of examples of weak categories.
In particular the many-to-one nerve functor sends strict $\o$-categories to multitopic categories.


The need to have a good description of higher many-to-one shapes was already clear at the conference {\em n-categories:  Foundations and Applications}
at IMA in Minneapolis, in June 2004.
Now there is (at least) seven essentially different definitions that are attempting to describe shapes of indets of many-to-one computads or some supposedly equivalent entities. These definitions differ a lot in spirit and it is by far not clear that they are all equivalent.  It seem that it is too early to call which one is better then the others and I think that all of them contribute to better understanding the concept they try to capture. So I will content myself by just listing them divided into
four groups.

1. There are three kinds of opetopes [BD], [L], [KJBM] that describe the set of shapes of many-to-one indets without an attempt to make it into a category.
The second and third kind of opetopes are proved in [KJBM] to be equivalent.

2. There are four categories describing the shapes of many-to-one indets: the category of multitopes, c.f. \cite{HMP}, the category of dendrotopes, c.f. \cite{Palm}, the category of opetopes \cite{Cheng} and the category of ordered face structures, c.f. \cite{Z2}. The main purpose of this paper is to show that the categories  first and last are equivalent.

3. The set of shapes of the, so called, pasting diagrams\footnote{$n$-pasting diagram are nothing but the shapes of domains of $n+1$ indets.} is described in \cite{HMP} as pasting diagrams and in \cite{Z2} as normal ordered face structures.

4. The category of all the shapes of many-to-one cells is the category $\ofs$ described in  \cite{Z2}.


The paper is organized as follows. In Section \ref{def_of_ofs} we recall the definition of an ordered face structures and two kinds of maps between them: monotone and local. In section \ref{GT-theory} we introduce the notion of a GT-theory that describes the abstract properties of $\ofs$. Sections \ref{ofsl} to \ref{third_adj} establish the main goal of the paper. Through a sequence of three adjunctions we establish that the category of multitopes and the category of principal face structures are equivalent. The remaining three sections are exhibiting the properties of $\ofs$. In Section \ref{shape} we define the category of shapes. In Section \ref{pra monads} we recall the relevant part of the work of T. Leinster and M. Weber in a way that is suitable for our context. Finally, in Section \ref{setup} we describe how the category $\ofs$ can generate all the ingredients involved in the definition of the many-to-one nerve construction for strict $\o$ categories.

As this paper is a sequel to \cite{Z2} we adopt here the notions
and notation introduced there.  This includes that we shall denote
the compositions of morphisms both ways, i.e. the composition of
two morphisms $X\stackrel{f}{\lra}Y\stackrel{g}{\lra}Z$ may be
denoted as either $g\circ f$ or $f;g$. But in any case we will
write which way the composition is meant.

I would like to thank the anonymous referee for comments that encouraged me to simplify the exposition
and to give a comprehensive explanation of the role the category of ordered face structures $\ofs$.
I also want to thank J. Kock for bringing \cite{Weber} to my attention.

The diagrams for this paper were prepared with a help of {\em catmac} of Michael Barr.

\section{Ordered face structures in a nutshell}\label{def_of_ofs}
This section is a quick introduction to ordered face structures.
For more see \cite{Z2}.

A {\em hypergraph}\index{hypergraph} $S$ is
\begin{enumerate}
  \item a family $\{ S_k\}_{k\in \o}$ of finite sets of faces;
  only finitely many among these sets are
  non-empty;
  \item a family of
functions $\{ \gamma : S_{k+1}\dsum \b1_{S_k}\ra S_k \}_{k \in\o
}$; where $\b1_{S_{k}}=\{ 1_u : u\in S_{k}\}$ is the set of {\em empty
faces} of dimension $k$; the face $1_u$ is the empty
$(k+1)$-dimensional face on a non-empty face $u$ of dimension $k$.
  \item a family of total relations $\{ \delta :
S_{k+1}\dsum \b1_{S_k}\ra S_k\dsum \b1_{S_{k-1}}\}_{k\in\o}$; for
$a\in S_{k+1}$ we denote $\delta(a)=\{ x\in S_k\dsum \b1_{S_{k-1}}
\, :\, (a,x)\in\delta^S_{k} \}$; $\delta(a)$ is either singleton
or it is non-empty subset of $S_k$\footnote{In other words $\delta(a)$ is either
equal to $\{ 1_x \}$ for some $x\in S_{k-1}$ or it
is a non-empty subset of $S_k$.}. Moreover $\delta :
S_{1}\dsum \b1_{S_0}\ra S_0$ is a function ($S_{-1}=\emptyset$).
We put $\dot{\delta}(a)=\delta(a)\cap S$.
\end{enumerate}

A {\em morphism of hypergraphs}\index{morphism! of
hypergraphs}\index{hypergraph!morphism} $f:S\lra T$ is a family of
functions $f_k : S_k \lra T_k$ that preserves $\gamma$ and
$\delta$ i.e., for $k \in\o$, $\gamma\circ f_{k+1}=f_k\circ\gamma$
and for $a\in S_{k+1}$ the restriction of $f_k$ to $\delta(a)$:
$f_a:\delta(a)\lra\delta(f(a))$ is a bijection (if $\delta(a)=1_u$
we mean by that $\delta(f(a))=1_{f(u)}$).

{\em Notation and conventions.}~If $a\in S_k$ we treat $\gamma(a)$
sometimes as an element of $S_{k-1}$ and sometimes as a subset $\{
\gamma(a) \}$ of $S_{k-1}$. Similarly $\delta(a)$ is treated
sometimes as a set of faces or as a single face if this set of
faces is a singleton. In particular, we say that a face $a$ is a
{\em loop} if $\gamma(a)=\delta(a)$ and by this we mean rather $\{
\gamma(a) \}=\delta(a)$. If $X$ is a set of faces in $S$ then by
$X^{-\lambda}$ we denote the set of faces in $X$ that are not
loops; $\dot{\delta}^{-\lambda}(a)=\dot{\delta}(a)\cap
S^{-\lambda}$. The set of {\em internal faces} of $a$ is
$\iota(a)=\gamma\dot{\delta}^{-\lambda}(a)\cap\delta\dot{\delta}^{-\lambda}(a)$.
The set $\theta(a)=\delta(a)\cup\gamma(a)$
($\dot{\theta}(a)=\dot{\delta}(a)\cup\gamma(a)$) is the sets of
(non-empty) faces of codimension $1$  in $a$.

On each set $S_k$ we introduce two binary relations $<^-$ and
$<^+$. On $S_0$ the relation $<^-$ is empty. If $k>0$, $<^-$ is
the transitive closure of the relation $\lhd^-$ on $S_k$, such
that $a \lhd^- b$ iff $\gamma(a)\in\delta(b)$.  We write
$a\perp^{S_k,-}b$ if either $a <^{S_k,-} b$ or $b <^{S_k,-} a$.
$<^+$ is the transitive closure of the relation  $\lhd^+$ on
$S_k$, such that $a \lhd^+ b$ iff $a\neq b$ and there is
$\alpha\in S_{k+1}^{-\lambda}$, such that $a\in \delta (\alpha)$
and $\gamma(\alpha)=b$.  We write $a\perp^{S_k,+}b$ if either $a
<^{S_k,+} b$ or $b <^{S_k,+} a$.

Let $A,B\subseteq S_k\cup \b1_{S_{k-1}}$. We set that $A$ is
1-{\em equal} $B$, notation $A\equiv_1B$, iff $A\cup
1_{\theta(A\cap S)} = B\cup 1_{\theta(B\cap S)}$.

An {\em ordered face structure}\index{face structure!ordered}
$(S,<^{S_k,\sim})_{k\in\o}$ (also denoted $S$) is a hypergraph $S$
together with a family of $\{<^{S_k,\sim}\}_{k\in\o}$ of binary
relations ($<^{S_k,\sim}$ is a relation on $S_k$),   if
$S_0\neq\emptyset$ and
\begin{enumerate}
 \item {\em Globularity:}\index{globularity}  for  $a\in S_{\geq 2}$:
$ \gamma\gamma(a)
=\gamma\delta(a)-\delta\dot{\delta}^{-\lambda}(a)$,
  $\delta\gamma(a) \equiv_1\delta\delta(a)-\gamma\dot{\delta}^{-\lambda}(a)$
  and for any $x\in S$:
  $\delta(1_x)=x=\gamma(1_x)$.
  \item {\em Local discreteness}:\index{local discreteness} if
  $x,y\in\delta(a)$ then $x\not\perp^+y$.
  \item {\em Strictness}:\index{strictness} for $k\in\o$, the
  relations $<^+$ and $<^\sim$ are  strict orders\footnote{By {\em strict order}
  we mean an irreflexive and transitive relation.} on $S_k$;
  $<^+$ on $S_0$ is linear.
  \item {\em Disjointness}:
  $\perp^\sim\cap \perp^+=\emptyset$, and
  for any $a,b\in S_k$:
  if $a<^\sim b$ then   $a<^- b$
 moreover if  $\theta(a)\cap\theta(b)= \emptyset$ then
  $a<^\sim b$ iff $a<^- b$.

  \item {\em Pencil linearity}: for any $a,b\in S_{\geq 1}$, $a\neq b$,
  \[ \mathrm{if}\;\; \dot{\theta}(a)\cap\dot{\theta}(b)\neq \emptyset\;\;\mathrm{then}\;\;
  \mathrm{either}\;\; a\perp^\sim b\;\; \mathrm{or}\;\; a\perp^+ b \]
  for any $a\in S_{\geq 2}$such that $\delta(a)\in\b1_S$, $b\in S_{\geq 2}$,
  \[ \mathrm{if}\;\; \gamma\gamma(a)\in\iota(b)\;\;\mathrm{then}\;\;
  \mathrm{either}\;\; a<^\sim b\;\; \mathrm{or}\;\; a<^+ b \]
  \item {\em Loop-filling}:
  \index{loop-filling} $S^\lambda\subseteq \gamma(S^{-\lambda})$ (where $S^\lambda$ is the
  set of loops in $S$ and $S^{-\lambda}=S-S^\lambda$).
\end{enumerate}

The {\em monotone morphism} of ordered face structures $f:S\lra T$
is a hypergraph morphism that preserves the order $<^\sim$. The
category of ordered face structures and monotone maps, is denoted by $\ofs$.

The {\em size of an ordered face structure} $S$ is the sequence
natural numbers $size(S)=\{ | S_n - \delta (S_{n+1}^{-\lambda})|
\}_{n\in\o}$, with almost all being equal $0$. We have an order
$<$ on such sequences, so that $\{ x_n \}_{n\in\o} < \{ y_n
\}_{n\in\o}$ iff there is $k\in\o$ such that $x_k< y_k$ and for
all $l>k$, $x_l = y_l$. This order is well founded and many facts
about ordered face structures can be proven by induction on the
size.  $S$ is {\em principal} iff $size(S)_l\leq1$, for $l\in\o$.
By $\pfs$ we denote full subcategory of $\ofs$ whose objects are
principal ordered face structures. In \cite{Z2} it was shown that
either an ordered face structure is principal or there is a cut
$\check{a}$ of $S$ that is defining a proper decomposition of $S$
into two ordered face structures $S^{\da \check{a}}$ and $S^{\ua
\check{a}}$ of smaller size than $S$ such that their $k$-tensor
product $S^{\da \check{a}}\otimes_kS^{\ua \check{a}}$ is
isomorphic to $S$, where $k$ is the dimension of the cut. By
$Sd(S)$ we denote the set of cuts of $S$ defining proper
decompositions of $S$.

The relation $<^\sim$  induces a binary relation
$(\dot{\delta}(a),<^{\sim}_a) $ for each $a\in S_{>0}$ (where
$<^{\sim}_a$ is the restriction of $<^\sim$ to the set
$\dot{\delta}(a)$). The {\em local morphism} of ordered face
structures $f:S\lra T$ is a hypergraph morphism that is a local
isomorphism  i.e. for $a\in S_{>1}$ the restricted map
$f_a:(\dot{\delta}(a),<^\sim_a)\lra(\dot{\delta}(f(a)),<^\sim_{f(a)})$
is an order isomorphism, where $f_a$ is the restriction of $f$ to
$\dot{\delta}(a)$. The category of ordered face structures and
local morphisms is denoted by $\ofsl$.

\section{Graded tensor theories}\label{GT-theory}

If we denote by $\ofs_n$ the full subcategory of  $\ofs$ containing object of dimension at most $n$
then for ordered face structures $S$, $S'$ we have operations of  the $k$-th domain $\bd^{(k)}S$, the $k$-th codomain
$\bc^{(k)}S$ and  $k$-tensor $S\otimes_kS'$ (whenever $\bc^{(k)}S=\bd^{(k)}S'$),  see \cite{Z2}. However
the operations $\bd^{(k)}S$ and $\bc^{(k)}S$ are not functorial with respect to the monotone morphisms
and $\ofs$ is not a monoidal globular category  in the sense of Batanin, see \cite{Batanin},
contrary to a public statement I have made.  On the other hand the category of
ordered face structures and inner maps, defined in Section \ref{setup}  is a monoidal globular category.

To explain the essential abstract structure of the category $\ofs$ I introduce
the notion of a graded tensor category. $\bN$ is the poset natural numbers.

{\em Graded tensor theory} $\cC$ (GT-theory for short) is a category $\cC$ equipped with
\begin{enumerate}
  \item a {\em dimension} functor $dim: \cC \ra \bN$;
  $\cC_k$ is the full subcategory of $\cC$ whose objects have dimension at most $k$;
  \item the objects of $\cC$ are rigid (i.e. no non-trivial automorphisms\footnote{In particular if two objects are isomorphic in $\cC$
   the isomorphism is always unique. This allow us treating isomorphic objects as equal and morphism with isomorphic domains and codomains as parallel.});
  \item for any object $S$ of $\cC$, such that $k\leq n=dim(S)$ there are
  domain and codomain morphisms
  \begin{center} \xext=800 \yext=450
\begin{picture}(\xext,\yext)(\xoff,\yoff)
\settriparms[-1`-1`0;400]
 \putAtriangle(0,50)[S`\bd^{(k)}S`\bc^{(k)}S;\bd^{(k)}_S`\bc^{(k)}_S`]
\end{picture}
\end{center}
 such that $dim(\bd^{(k)}S)\leq k= dim(\bc^{(k)}S)$ and,  for $k< l\leq n$,
 \[ \bd^{(l)}_S\circ \bc^{(k)}_{\bd^{(l)}S} =\bc^{(l)}_S\circ \bc^{(k)}_{\bc^{(l)}S} =\bc^{(k)}_S \]
  \[ \bc^{(l)}_S\circ \bd^{(k)}_{\bc^{(l)}S} =\bd^{(l)}_S\circ \bd^{(k)}_{\bd^{(l)}S} =\bd^{(k)}_S \]
in particular the diagram
  \begin{center}
\xext=1000 \yext=950 \adjust[`I;I`;I`;`I]
\begin{picture}(\xext,\yext)(\xoff,\yoff)
 \settriparms[-1`-1`0;500]
 \putAtriangle(0,500)[S`\bd^{(l)}S`\bc^{(l)}S;\bd^{(l)}_S`\bc^{(l)}_S`]
 \setsqparms[0`-1`-1`0;1000`500]
\putsquare(0,0)[\phantom{\bd^{(l)}S}`\phantom{\bc^{(l)}S}`\bd^{(k)}S`\bc^{(k)}S; `\bd^{(k)}_{\bd^{(l)}S}`\bc^{(k)}_{\bc^{(l)}S}`]

 \put(200,50){\vector(2,1){700}}
 \put(450,225){\vector(-2,1){350}}
 \put(800,50){\line(-2,1){250}}
 \put(640,380){\makebox(100,100){$\bd^{(k)}_{\bc^{(l)}S}$}}
 \put(240,380){\makebox(100,100){$\bc^{(k)}_{\bd^{(l)}S}$}}
\end{picture}
\end{center}
commutes. If $k\geq n$, we put   $\bd^{(k)}S$ and  $\bc^{(k)}S$ to be identities on $S$.
  \item For $k< n$, the category $\cC_n \times_k\cC_n$ is the category whose objects consists of three objects, two in $\cC_n$ and one in $\cC_k$
  and two maps as follows
  \begin{center} \xext=800 \yext=450
\begin{picture}(\xext,\yext)(\xoff,\yoff)
\settriparms[0`-1`-1;400]
 \putVtriangle(0,50)[S`S'`\bd^{(k)}S\cong\bc^{(k)}S';`\bc^{(k)}_S`\bd^{(k)}_{S'}]
\end{picture}
\end{center}
Clearly the objects of $\cC_n \times_k\cC_n$ can be thought of as pairs $(S,S')$ of objects of $\cC_n$ satisfying an obvious compatibility condition.
The morphisms in $\cC_n \times_k\cC_n$ are triples of morphisms in $\cC$ commuting with the morphisms $\bc^{(k)}$ and $\bd^{(k)}$.
We have three functors
\[ \otimes_k,\pi_1,\pi_2 : \cC_n \times_k\cC_n \lra \cC_n \]
where $\pi_1$ and $\pi_2$ are the obvious projections and two natural transformations
 \begin{center}
\xext=1000 \yext=100
\begin{picture}(\xext,\yext)(\xoff,\yoff)
\putmorphism(0,20)(1,0)[\pi_1`\otimes_k `\kappa^1]{500}{1}a
\putmorphism(500,20)(1,0)[\phantom{\otimes_k}`\pi_2 `\kappa^2]{500}{-1}a
\end{picture}
\end{center}
so that the squares
\begin{center}
\xext=800 \yext=600
\begin{picture}(\xext,\yext)(\xoff,\yoff)
 \setsqparms[1`-1`-1`1;800`500]
 \putsquare(0,50)[S`S\otimes_k S'`\bc^{(k)}S`S';\kappa^1_S`\bc^{(k)}_{S}`\kappa^2_{S'}`\bd^{(k)}_{S'}]
 \end{picture}
\end{center}
commute, for any $(S,S')$ in  $\cC_n \times_k\cC_n$. Such squares, or squares isomorphic to them, are called {\em $k$-tensor squares} or simply {\em tensor squares}.

\item These data are related by the existence of some isomorphisms between objects and morphisms.
But, as $\cC$ is rigid, these isomorphisms are necessarily unique and hence there are no coherence conditions since any diagram of isomorphisms in $\cC$ commutes. Thus we will put no names on those isomorphisms. $R,R',S,S',S''$ are assumed to be ordered face structures.
    \begin{enumerate}
      \item {\em Domains and codomains of compositions.} For $k>l$, there are isomorphism making the triangles
      \begin{center}
\xext=3000 \yext=550
\begin{picture}(\xext,\yext)(\xoff,\yoff)
 \settriparms[-1`-1`1;450]
 \putAtriangle(250,50)[S\otimes_l S'`\bd^{(k)}S\otimes_l \bd^{(k)}S'`\bd^{(k)}(S\otimes_l S');
 \bd^{(k)}_S\otimes_l \bd^{(k)}_{S'}`\bd^{(k)}_{S\otimes_l S'}`\cong]
  \putAtriangle(1900,50)[S\otimes_l S'`\bc^{(k)}S\otimes_l \bc^{(k)}S'`\bc^{(k)}(S\otimes_l S');
 \bc^{(k)}_S\otimes_l \bc^{(k)}_{S'}`\bc^{(k)}_{S\otimes_l S'}`\cong]
\end{picture}
\end{center}
commute. For $k\leq l$, there are isomorphism making the triangles
\begin{center}
\xext=3000 \yext=550
\begin{picture}(\xext,\yext)(\xoff,\yoff)
 \settriparms[-1`-1`1;450]
 \putAtriangle(250,50)[S\otimes_l S'`\bd^{(k)}S`\bd^{(k)}(S\otimes_l S');
 \kappa^1_S\circ \bd^{(k)}_S`\bd^{(k)}_{S\otimes_l S'}`\cong]
  \putAtriangle(1900,50)[S\otimes_l S'`\bc^{(k)}S'`\bc^{(k)}(S\otimes_l S');
 \kappa^2_{S'}\circ \bc^{(k)}_S`\bc^{(k)}_{S\otimes_l S'}`\cong]
\end{picture}
\end{center}
      commute.
      \item {\em Units.} The following diagrams
      \begin{center} \xext=2000 \yext=580
\begin{picture}(\xext,\yext)(\xoff,\yoff)
 \setsqparms[1`-1`-1`1;700`450]
 \putsquare(0,50)[X`X`\bd^{(k)}X`\bd^{(k)}X;1_T`\bd^{(k)}_T`\bd^{(k)}_T`1_{\bd^{(k)}X}]
  \putsquare(1300,50)[X`X`\bc^{(k)}X`\bc^{(k)}X;1_T`\bc^{(k)}_T`\bc^{(k)}_T`1_{\bc^{(k)}X}]
\end{picture}
\end{center}
are $k$-tensor squares.
      \item {\em Associativity.} Whenever any of the expressions make sense  we have an isomorphism
      \[ S\otimes_k(S'\otimes_k S'')\cong (S\otimes_k S')\otimes_kS''\]
      \item {\em Middle exchange.} For $k<l$, if we have a diagram
      \begin{center}
\xext=1200 \yext=800
\begin{picture}(\xext,\yext)(\xoff,\yoff)
   \putmorphism(0,400)(1,0)[\bc^{(l)}R`\bc^{(k)}R`\bc^{(k)}]{600}{-1}a
   \putmorphism(600,400)(1,0)[\phantom{\bc^{(k)}R}`\bc^{(l)}S`\bd^{(k)}]{600}{1}a

   \putmorphism(0,400)(0,1)[\phantom{\bc^{(l)}R}`R'`\bd^{(l)}]{400}{1}l
   \putmorphism(0,800)(0,1)[R`\phantom{\bc^{(l)}R}`\bc^{(l)}]{400}{-1}l

   \putmorphism(1200,400)(0,1)[\phantom{\bc^{(l)}R'}`S'`\bd^{(l)}]{400}{1}r
   \putmorphism(1200,800)(0,1)[S`\phantom{\bc^{(l)}R}`\bc^{(l)}]{400}{-1}r
\end{picture}
\end{center}
      then the two objects we can form out of it
      \[ (R\otimes_l R')\otimes_k (S\otimes_l S') \cong (R\otimes_k S)\otimes_l (R'\otimes_k S') \]
      are isomorphic.
    \end{enumerate}
\end{enumerate}

A {\em model} of a GT-theory $\cC$  in a category $\cA$, or $\cA$-model of $\cC$ for short,
is a functor from $\cC^{op}$ to $\cA$ which sends tensor squares to pullbacks and
the distinguished isomorphisms to the canonical isomorphisms. By $Mod_\otimes(\cC^{op},\cA)$ w denote
the category of $\cA$-models of $\cC$ and natural transformations.
A model $\cG_\cC:\cC^{op}\ra \cA$ is {\em generic} iff
for any other model $M:\cC^{op}\ra \cB$ there is a unique (up to an iso) functor $\overline{M}:\cA\ra \cB$ making the triangle
\begin{center}
\xext=500 \yext=550
\begin{picture}(\xext,\yext)(\xoff,\yoff)
 \settriparms[1`1`1;500]
 \putqtriangle(0,0)[\cC`\cA`\cB;\cG_\cC`M`\overline{M}]
\end{picture}
\end{center}
commute up to an isomorphism, and for any natural transformation $\sigma:M\ra M'$ between models there is a unique
 natural transformation $\overline{\sigma}:\overline{M}\ra \overline{M'}$ such that $\overline{\sigma}_{\cG_\cC} =\sigma$.
If identity functor on GT-theory $\cC$ is a model then it is the generic model of $\cC$ and $\cC$ is a called a {\em realized} GT-theory.

\vskip 2mm
{\em Examples.}
\begin{enumerate}
\item The category $\Delta_0$ of finite linear graphs is clearly a GT-theory. It contains one object $[0]$ of dimension $0$ and all
the other objects are of dimension $1$. The domain and codomain maps $\bd^{(0)}_{[n]},\bc^{(0)}_{[n]}:[0]\ra [n]$  send the only vertex
of $[0]$ to the first and last vertex in the in $[n]$, respectively. The $0$-tensor of any two
objects is defined and we have $[n]\otimes_0[m]=[n+m]$ with the obvious inclusions. In this case the tensor
squares are actual pushouts in $\Delta_0$, i.e. $\Delta_0$ is a realized GT-theory. Note that this tensor operation does not make
$\Delta_0$ a monoidal category as $+$ is not 'sufficiently functorial'.
\item  The category $s\oG$ of simple $\o$-graphs (or globular cardinals)  is also a realized GT-theory. If we look
at the objects of $s\oG$ as one-to-one pasting diagrams then the domain and the codomain operations are the
pasting diagrams of the $k$-th domain and  the $k$-th codomain of this diagram.
\item The whole GT-theory structure of the category $\fs$ of positive face structures is described in \cite{Z1}
and even in this case the tensor squares are pushouts, i.e. $\fs$ is a realized GT-theory, as well.
\item The GT-theory structure of the category $\ofs$ of ordered face structures is described in \cite{Z2}
but in this case the tensor squares are not pushouts in general. This is because only part of the order $<^\sim$ in the tensor
is determined by the components, see \cite{Z2} details. The embedding functor $\cG_\ofs:\ofs^{op}\ra\ofsl^{op}$ is the generic model of $\ofs$.
\item In Section \ref{setup} we shall define still another GT-theory $\ofsm$ of ordered face structures and monotone $\o$-maps which
has a non-identity generic model $\cG_{\ofsm}:\ofsm\ra\ofso$.
\end{enumerate}

We have the following
\begin{proposition}\label{GTmodel}Let $\cC$ be an essentially small GT-theory and $\cA$ be a locally small category, and
$M:\cC\lra \cA$ be a functor such that its dual $M^{op}:\cC^{op}\lra \cA^{op}$ is an $\cA^{op}$-model of $\cC$.
Then $M$ induces a functor $\bar{M}: \cA\lra \oC$.
\end{proposition}

{\it Proof.}~  I will give the construction of $\bar{M}$ only.  Let $A$ be an object of $\cA$. The $n$-cells of the $\o$-category $\bar{M}(A)$ are given by
\[ \bar{M}(A)_n=\coprod_{S} \cA(M(S),A) \]
where the coproduct is taken over all (up to isomorphism\footnote{In practice we think that we take the coproduct over all objects of $\cC$ of dimension at most $n$ and we identify two maps $a:M(S)\ra A$ with  $a':M(S')\ra A$ if there is a (necessarily unique) isomorphism $i:S:\ra S'$ such that $a'\circ M(i)=a$.}) objects of $\cC$ of dimension at most $n$. If $S$ has dimension lower than $n$ than the morphism $M(S)\ra A$ is considered as the identity at dimension $n$ of a lower dimension cell. In particular, the identity operations in the $\o$-category $\bar{M}$ are inclusions. The $k$-th domain and $k$-th codomain  operations \mbox{$d^{(k)},c^{(k)}: \bar{M}(A)_n\lra\bar{M}(A)_k$}  are defined by composition as follows. For \mbox{$a:M(S)\ra A\in \bar{M}(A)_n$} we put
\[ d^{(k)}(a)= a\circ M(\bd^{(k)}_S),\hskip 5mm c^{(k)}(a)= a\circ M(\bc^{(k)}_S). \]
If $b:S'\ra A\in \bar{M}(A)_n$ so that $c^{(k)}(a)=d^{(k)}(b)$ the outer square in the following diagram
\begin{center} \xext=1500 \yext=1150
\begin{picture}(\xext,\yext)(\xoff,\yoff)
 \setsqparms[1`-1`-1`1;900`600]
 \putsquare(0,50)[M(S)`M(S\otimes_kS')`M(\bc^{(k)}S)`M(S');`M(\bc^{(k)}_S)``M(\bd^{(k)}_{S'})]
  \put(900,750){\makebox(100,100){$a{\bf ;}_kb$}}
  \put(500,800){\makebox(100,100){$a$}}
  \put(1200,350){\makebox(100,100){$b$}}
  \put(1500,975){\makebox(100,100){$A$}}
 \put(950,700){\vector(2,1){550}}
 \put(100,700){\vector(4,1){1400}}
 \put(1000,100){\vector(2,3){550}}
 \end{picture}
\end{center}
commutes. Since $M$ is a model of $\cC$ the inner square is a pushout
and we have a morphism $a{\bf ;}_kb:M(S\otimes S')\ra A$ making the remaining two triangles commute.
$~\Box$

\vskip 2mm
{\em Remark.} It is convenient to thing about $\ofs$ in terms of the abstract
notion of a graded tensor category. In fact more than the above Lemma can be stated for an
abstract graded tensor category not only $\ofs$. But I think that the general theory of GT-theories should
wait until more non-trivial GT-theories that are not a realized GT-theories are found.

\section{The category $\ofsl$}\label{ofsl}

The following Lemma subsume some properties of the category of
ordered face structures and local maps $\ofsl$ that are
essentially in \cite{Z2}.

\begin{lemma}
Let $f:S\ra T$ and $g: P\ra T$ be morphisms in $\ofsl^{op}$, with
$P$ being a principal ordered face structure of dimension $n$,
$P_n=\{ \bm_P\}$, $a\in S_n$. If $f(a)=g(\bm_P)$ then there is a
unique map $\bar{g}:P\ra S$ such that $\bar{g}(\bm_P)=a$ and hence
$f\circ \bar{g}=g$. In particular, any principal ordered face
structure is projective in $\ofsl^{op}$.
\end{lemma}

{\it Proof.}~  The first statement follows from Lemmas 11.1 and
11.2 from \cite{Z2}. To see that this imply that principal ordered
face structures are projective in $\ofsl^{op}$ it is enough to
note that the local maps in $\ofsl$ are epi iff they are onto.
$~\Box$

Let $S$ be an ordered face structure. We have an obvious
projection functor
\[ \Sigma^S : \pfs \da S \lra  \pfs \lra \ofsl \]
such that
\[ \Sigma^S (f:P\ra S) = P \]
and the {\em principal cocone over} $S$
\[ \sigma^S : \Sigma^S\lra S\]
such that
\[ \sigma^S_{(f:P\ra S)} = f : \Sigma^S (f:P\ra S) = P\lra S \]

We have

\begin{lemma}\label{speclim}
The cocone $\sigma^S: \Sigma^S \stackrel{\cdot}{\lra} S$ is a
colimiting cocone in $\ofsl$.
\end{lemma}

{\it Proof.}~  Simple check. $~\Box$

The colimits in $\ofsl$, as described above, are called {\em
principal colimits} (and when considered in $\ofsl^{op}$ are
called {\em principal limits}). Note that we are not saying that
the above cocone is a colimit in $\ofs$ i.e. in the category of
ordered face structures and monotone maps. For example if $S$ is
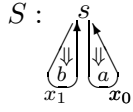
\begin{center} \xext=1500 \yext=380
\begin{picture}(\xext,\yext)(\xoff,\yoff)
  \put(440,300){$S: $}
  \put(790,0){$^{x_0}$}
  \put(630,130){\oval(100,100)[b]}
  \put(580,130){\line(1,2){85}}
  \put(680,130){\vector(0,1){180}}
  \put(620,130){$^\Da$}
    \put(610,65){$^{b}$}
     \put(565,0){$^{x_1}$}
     \put(680,320){$s$}

  \put(770,130){\oval(100,100)[b]}
   \put(720,130){\line(0,1){180}}
  \put(820,130){\vector(-1,2){85}}
   \put(730,130){$^\Da$}
   \put(750,65){$^{a}$}
   \put(790,0){$^{x_0}$}
\end{picture}
\end{center}
then clearly the principal cocone over $S$ is not a colimiting
cocone in $\ofs$ as it does not determine the order $<^\sim$ between faces
$x_1$ and $x_0$. In fact the ordered face structures $S$ for which
the cocone $\sigma^S: \Sigma^S \stackrel{\cdot}{\lra} S$ is a
colimiting cocone in $\ofs$ have several good properties that are
going to be studied elsewhere.

The following Lemma states some properties of principal cocones
over tensors.

\begin{lemma} \label{tensor-principal cocone}
Let $S$ and $S'$ be ordered face structures such that
$\bc^{(k)}(S)=\bd^{(k)}(S')$, $P$ a principal ordered face
structure, and $f:P\lra S\otimes_kS'$ a map in $\ofsl$. Then
\begin{enumerate}
  \item either $f$ factorizes (uniquely) via $\kappa^1_S$ or $f$ factorizes (not necessarily uniquely) via
  $\kappa^2_{S'}$;
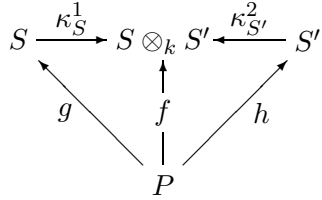
  \begin{center} \xext=2400 \yext=600
\begin{picture}(\xext,\yext)(\xoff,\yoff)
\settripairparms[1`-1`-1`-1`-1;500]
\putVtrianglepair(0,0)[S`S\otimes_kS'`S'`P;\kappa^1_S`\kappa^2_{S'}`g`f`h]
 \putmorphism(600,500)(1,0)[\phantom{S'}`S'`\kappa^2_{S'}]{400}{-1}a
\end{picture}
\end{center}
  \item if there are both $g$ and $h$ factorizations of $f$ then
  there is a factorization $l$ making the diagram
    \begin{center} \xext=2400 \yext=600
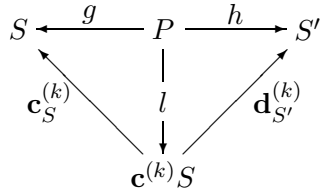
\begin{picture}(\xext,\yext)(\xoff,\yoff)
\settripairparms[-1`-1`-1`1`-1;500]
\putVtrianglepair(0,0)[S`P`S'`\bc^{(k)}S;g`h`\bc^{(k)}_S`l`\bd^{(k)}_{S'}]
 \putmorphism(500,500)(1,0)[\phantom{S'}`S'`h]{500}{1}a
\end{picture}
\end{center}
   commute;
  \item finally, if there are two factorizations $h$ and $h'$ of $f$ via $\kappa^2_{S'}$
  then there a factorization $g$ via $\kappa^1_S$.
\end{enumerate}
\end{lemma}

{\it Proof.}~  This easily follows from the explicite description
of the tensors in  \cite{Z2}.
 $~\Box$

\section{Simple adjunction}
The proof of the main theorem proceeds by establishing three adjoint equivalences.

By Lemma 11.1 of [Z2], the inclusion functor $\bi : \pfs \lra
\ofs_{loc}$ is full and faithful. It induces the adjunction
\begin{center} \xext=1500 \yext=250
\begin{picture}(\xext,\yext)(\xoff,\yoff)
\putmorphism(0,150)(1,0)[\phantom{Set^{{\pfs}^{op}}}`
\phantom{Set^{{(\ofs)}^{op}}}`Ran_\bi]{1500}{1}a
\putmorphism(0,50)(1,0)[\phantom{Set^{\pfs^{op}}}`
\phantom{Set^{{(\ofs)}^{op}}}`\bi^*]{1500}{-1}b
\putmorphism(0,100)(1,0)[Set^{{\pfs}^{op}}`Set^{\ofs_{loc}^{op}}`]{1500}{0}b
\end{picture}
\end{center}
where $\bi^*$ is the functor of composing with $\bi$ and $Ran_\bi$
is the right Kan extension along $\bi$. Recall that for $F$ in
$Set^{\pfs^{op}}$, $S$ in $\ofs$, it is defined as the following
limit
\[ (Ran_\bi F)(S) = Lim F\circ (\Sigma^{S})^{op} \]
where $(\Sigma^{S})^{op}:S\da \pfs^{op}\ra\pfs^{op}$.  Clearly
$Ran_\bi F$ preserves principal limits. As $\bi$ is full and
faithful the right Kan extension $Ran_\bi (F)$ is an extension,
i.e. the counit of this adjunction
\[ \varepsilon_F : (Ran_\bi\, F)\circ \bi \lra F \]
is an isomorphism. In particular, $Ran_\bi$ is full and faithful.
It is easy to see, that for $G$ in $Set^{\ofsl^{op}}$ the unit of
adjunction
\[  \eta_G : G \lra Ran_\bi (G\circ\bi) \]
is an isomorphism iff $G$ preserves principal limits. Thus we have proved

\begin{proposition}\label{presheaf_eq1}
The above adjunction restricts to the following equivalence of
categories
\begin{center} \xext=1500 \yext=250
\begin{picture}(\xext,\yext)(\xoff,\yoff)
\putmorphism(0,150)(1,0)[\phantom{Set^{\pfs^{op}}}`
\phantom{Mod_\otimes((\ofs)^{op},Set)}`Ran_\bi]{1500}{1}a
\putmorphism(0,50)(1,0)[\phantom{Set^{\pfs^{op}}}`
\phantom{Mod_\otimes((\ofs)^{op},Set)}`\bi^*]{1500}{-1}b
\putmorphism(0,100)(1,0)[Set^{\pfs^{op}}`pLim(\ofsl^{op},Set)`]{1500}{0}b
\end{picture}
\end{center}
where $pLim(\ofsl^{op},Set)$ is the category of principal limits preserving functors and natural transformation.
$~~\Box$
\end{proposition}

\section{Tensor squares vs principal limits}

We shall define an adjunction
\begin{center} \xext=1200 \yext=220
\begin{picture}(\xext,\yext)(\xoff,\yoff)
\putmorphism(0,100)(1,0)[Set^{\ofsl^{op}}`Mod_\otimes(\ofs^{op},Set)`]{1200}{0}a
\putmorphism(0,150)(1,0)[\phantom{Set^{\ofsl^{op}}}`\phantom{Mod_\otimes(\ofs^{op},Set)}`\be]{1200}{-1}a
\putmorphism(0,50)(1,0)[\phantom{Set^{\ofsl^{op}}}`\phantom{Mod_\otimes(\ofs^{op},Set)}`\cL]{1200}{1}b
\end{picture}
\end{center}
The functor $\be$ is sending $Set$-models of $\ofs$ to the presheaves on $\ofsl$ along the generic model
$\cG_\ofs$ and natural transformations to the same natural transformations. Thus $\be$ can be thought of
as an embedding that is extending the models of $\ofs$ by defining them on all the local maps in $\ofsl$.
We tend to omit $\be$ in formulas writing for example the unit of this
adjunction as $\eta_F : F\ra \cL(F)$,  understanding that it is a
morphism in $Set^{\ofsl^{op}}$ rather than in
$Mod_\otimes(\ofs^{op},Set)$.

Let $F: \ofsl^{op}\ra Set$ be a presheaf. Then
\[ \cL(F) = Lim(F\circ (\Sigma^{(-)})^{op}) : \ofsl^{op}\lra Set  \]
i.e. $\cL(F)(S)$ is the limit of the following functor
\begin{center} \xext=1600 \yext=100
\begin{picture}(\xext,\yext)(\xoff,\yoff)
\putmorphism(0,0)(1,0)[(\pfs\da S)^{op}`\ofsl^{op}`(\Sigma^S)^{op}]{1000}{1}a
\putmorphism(1000,0)(1,0)[\phantom{\pfs^{op}}`Set`F]{800}{1}a
\end{picture}
\end{center}
with the limiting cone $\sigma^{F,S}: \cL(F)(S)\lra F\circ
(\Sigma^S)^{op}$. For a monotone map \mbox{$f:S\ra S'$} the function
$\cL(F)(f)$ is so defined that for any $h:P\ra S$ in
$(\pfs\da S)^{op}$ the triangle
\begin{center} \xext=900 \yext=550
\begin{picture}(\xext,\yext)(\xoff,\yoff)
 \settriparms[1`1`1;450]
 \putVtriangle(0,0)[\cL(F)(S)`\cL(F)(S')`F(P);\cL(F)(f)`\sigma^{F,S}_h`\sigma^{F,S'}_{f\circ h}]
\end{picture}
\end{center}
commutes. For a natural transformation $\tau : F\lra F'$ and $S\in
\ofs$, $\cL(\tau)_S$ is the unique map making the squares
\begin{center} \xext=800 \yext=600
\begin{picture}(\xext,\yext)(\xoff,\yoff)
 \setsqparms[1`1`1`1;800`500]
 \putsquare(0,0)[\cL(F)(S)`\cL(F')(S)`F(P)`F'(P);\cL(\tau)_S`\sigma^{F,S}_h`\sigma^{F',S}_h`\tau_P]
\end{picture}
\end{center}
commutes, for any $h:P\ra S\in \ofs$.

\begin{proposition} \label{adj1}
$\cL $ is well defined functor and $\cL\dashv \be$.
\end{proposition}

{\it Proof.}~The fact that $\cL(F)$ is a functor and $\cL(\tau)$
is a natural transformation is left to the reader. We shall verify
that $\cL(F)$ is a model of $\ofs$, i.e. sends tensor squares to pullbacks.

Let $S$ and $S'$ be ordered face structures such that
$\bc^{(k)}S=\bd^{(k)}S'$.  As $\cL$ is a functor it sends commuting
squares to commuting squares and hence we have a unique function
$\varphi$ making the diagram
\begin{center} \xext=2400 \yext=950
\begin{picture}(\xext,\yext)(\xoff,\yoff)
 \putmorphism(0,400)(1,0)[\cL(F)(S)`\cL(F)(S)\times_{\cL(F)(\bc^{(k)}S)}\cL(F)(S')`\pi^{F,S}]{1200}{-1}a
 \putmorphism(1200,400)(1,0)[\phantom{\cL(F)(S)\times_{\cL(F)(\bc^{(k)}S)}\cL(F)(S')}`\cL(F)(S')`\pi^{F,S'}]{1200}{1}a

\putmorphism(0,400)(3,-1)[\phantom{\cL(F)(S)}`\cL(F)(\bc^{(k)}S)`]{1200}{1}l
\putmorphism(0,350)(3,-1)[\phantom{\cL(F)(S)}`\phantom{\cL(F)(\bc^{(k)}S)}`\cL(F)(\bc^{(k)}_S)]{1200}{0}l

\putmorphism(2400,400)(-3,-1)[\phantom{\cL(F)(S)}`\phantom{\cL(F)(\bc^{(k)}S)}`]{1200}{1}r
\putmorphism(2400,350)(-3,-1)[\phantom{\cL(F)(S)}`\phantom{\cL(F)(\bc^{(k)}S)}`\cL(F)(\bd^{(k)}_{S'})]{1200}{0}r

\putmorphism(1200,900)(0,-1)[\cL(F)(S\otimes_kS')``\varphi]{500}{1}r

\putmorphism(1000,900)(-2,-1)[\phantom{\cL(F)(S\otimes_kS')}`\phantom{\cL(F)(S)}`\cL(F)(\kappa^1_S)]{1000}{1}l
\putmorphism(1400,900)(2,-1)[\phantom{\cL(F)(S\otimes_kS')}`\phantom{\cL(F)(S')}`\cL(F)(\kappa^2_{S'})]{1000}{1}r
\end{picture}
\end{center}
commute. We shall define a function
\[ \psi :\cL(F)(S)\times_{\cL(F)(\bc^{(k)}S)}\cL(F)(S')\lra \cL(F)(S\otimes_kS')\]
the inverse of $\varphi$, by defining a cone $\xi$ from
$\cL(F)(S)\times_{\cL(F)(\bc^{(k)}S)}\cL(F)(S')$ to the functor
\begin{center} \xext=2000 \yext=100
\begin{picture}(\xext,\yext)(\xoff,\yoff)
\putmorphism(0,0)(1,0)[(\pfs\da(S\otimes_kS'))^{op}`\ofsl^{op}`(\Sigma^{S\otimes_kS'})^{op}]{1200}{1}a
\putmorphism(1200,0)(1,0)[\phantom{\pfs^{op}}`Set`F]{800}{1}a
\end{picture}
\end{center}
Let $f:P\ra S\otimes_kS'$ be a map in $\ofsl$
\[\xi_f = \left\{ \begin{array}{ll}
             \pi^{F,S};\sigma^{F,S}_g &  \mbox{if $f=g;\kappa^1_S$ for some $g:P\ra S$,}  \\
             \pi^{F,S'};\sigma^{F,S'}_h &  \mbox{if $f=h;\kappa^2_{S'}$ for some $h:P\ra S'$,}
                                    \end{array}
                \right.\]
The Lemma \ref{tensor-principal cocone} guarantee that this
definition gives in fact a cone $\xi$ over
\mbox{$F\circ(\Sigma^{S\otimes_kS'})^{op}$}. Thus we have $\psi$ as in the
diagram
\begin{center} \xext=2400 \yext=760
\begin{picture}(\xext,\yext)(\xoff,\yoff)
 \putmorphism(300,670)(1,-4)[\cL(F)(S)\times_{\cL(F)(\bc^{(k)}S)}\cL(F)(S')`\cL(F)(S)`\pi^{F,S}]{150}{1}l
 \putmorphism(300,670)(2,-1)[\phantom{\cL(F)(S)\times_{\cL(F)(\bc^{(k)}S)}\cL(F)(S')}`\cL(F)(S')`]{750}{1}r
 \put(750,470){$\pi^{F,S'}$}

  \putmorphism(1050,320)(3,-1)[\phantom{\cL(F)(S')}``]{850}{1}r
  \put(1450,220){$\sigma^{F,S'}_h$}

 \putmorphism(2200,670)(0,-1)[\cL(F)(S\otimes_kS')`F(Q)=F\circ(\Sigma^{S\otimes_kS'})^{op}(f)`
 \sigma^{F,S\otimes_kS'}_f ]{600}{1}r

 \putmorphism(480,70)(1,0)[\phantom{\cL(F)(S)}`\phantom{F(Q)=F\circ(\Sigma^{S\otimes_kS'})^{op}(f)}`
 \sigma^{F,S}_g]{1700}{1}b

 \putmorphism(300,670)(1,0)[\phantom{\cL(F)(S)\times_{\cL(F)(\bc^{(k)}S)}\cL(F)(S')}`
  \phantom{\cL(F)(S\otimes_kS')}`\psi]{1900}{1}a
\end{picture}
\end{center}
As both $\varphi$ and $\psi$ are defined using universal
properties of limits, it is easy to see that they are mutually
inverse, i.e. $\cL(F)$ preserves special pullbacks.

As the image under $F$ of the principal cocone $\sigma^S$ from
$\Sigma^S$ to $S$ is a cone from $F(S)$ to $F\circ
(\Sigma^S)^{op}$ we have a unique  map $(\eta_F)_S:F(S)\lra
\cL(F)(S)$ making the triangles
\begin{center} \xext=900 \yext=520
\begin{picture}(\xext,\yext)(\xoff,\yoff)
\settriparms[1`1`1;400] \putVtriangle(0,0)[F(S)`\cL(F)(S)`F(P);
(\eta_F)_S`F(h)`\sigma^{F,S}_{h}]
\end{picture}
\end{center}
commute, for any $h:P\ra S\in \ofsl$. That defines the unit of
adjunction $\cL\dashv \be$. For any principal ordered face
structure $P$, the category $\pfs\da P$ has the terminal object
$1_P:P\ra P$. Thus any $P$-component of the unit of adjunction
$(\eta_F)_P:F(P)\lra \cL(F)(P)$ is an isomorphism.

The counit of adjunction $\varepsilon_G:\cL(G)\ra G$ is defined
using the fact that both $G$ and $\cL(G)$ are models of $\ofs$.
The map $(\varepsilon_G)_S:\cL(G)(S)\ra G(S)$ is defined by
induction on the size of $S$. If $S=P$ is a principal ordered face
structure then we put $(\varepsilon_G)_P=(\eta_F)_P^{-1}$. If $S$
is not principal than with $\check{a}\in Sd(S)_k$, and we have
$S=S^{\da \check{a}}\otimes_k S^{\da \check{a}}$ we put
 \[ (\varepsilon_G)_S= (\varepsilon_G)_{S^{\da \check{a}}}
 \times_{(\varepsilon_G)_{\bc^{(k)}(S^{\da \check{a}})}}  (\varepsilon_G)_{S^{\ua \check{a}}} \]
To verify the triangular equalities it is enough to show that the
triangles
\begin{center} \xext=2500 \yext=650
\begin{picture}(\xext,\yext)(\xoff,\yoff)
\settriparms[-1`1`1;500]
\putAtriangle(0,50)[\cL(\cL(F))(S)`\cL(F)(S)`\cL(F)(S);
\cL((\eta_F)_S)`(\varepsilon_{\cL(F)})_S)`1_{\cL(F)(S)}]
\putAtriangle(1600,50)[\cL(G)(S)`G(S)`G(S);
(\eta_G)_S`(\varepsilon_G)_S`1_{G(S)}]
\end{picture}
\end{center}
commute, for each ordered face structure $S$ separately. The
commutation of the left triangle can be shown using the fact that
all functors involved preserves principal limits and the
commutation of the right triangle can be shown by induction on the
size of $S$ using the fact that all the involved functors are
models of $\ofs$. The remaining details are left for the reader.
$~\Box$

\begin{proposition}\label{presheaf_eq2}
The above adjunction restricts to the following equivalence of
categories
\begin{center} \xext=1400 \yext=220
\begin{picture}(\xext,\yext)(\xoff,\yoff)
\putmorphism(0,100)(1,0)[pLim(\ofsl^{op},Set)`Mod_\otimes(\ofs^{op},Set)`]{1400}{0}a
\putmorphism(0,150)(1,0)[\phantom{pLim(\ofsl^{op},Set)}`\phantom{Mod_\otimes(\ofs^{op},Set)}`\be]{1400}{-1}a
\putmorphism(0,50)(1,0)[\phantom{pLim(\ofsl^{op},Set)}`\phantom{Mod_\otimes(\ofs^{op},Set)}`\cL]{1400}{1}b
\end{picture}
\end{center}
\end{proposition}
{\it Proof.}~ As $\be$ is full and faithful the counit of the adjunction $\be\dashv \cL$ is an isomorphism.
 From the description of the functor $\cL$ it is clear that, for any functor \mbox{$F:\ofsl^{op}\ra Set$},
  $\cL(F)$ preserves principal limits and that $\eta_F:F\ra \cL(F)$ is an isomorphism iff $F$
preserves principal limits.
$~\Box$

\section{Third adjunction}\label{third_adj}

Recall that in \cite{Z2} we have defined a functor  \mbox{$(-)^*:\ofs\ra\mComp$}
associating to any ordered face structure $S$ the many-to-one computad $S^*$ generated by $S$.
The $n$-cells of $S^*$ are (equivalence classes of) local maps $a:R\ra S$ from ordered face structures $R$ of dimension
at most $n$. The $k$-domain and the $k$-codomain of $a$ are $d^{(k)}(a)=a\circ \bd^{(k)}_R$ and $c^{(k)}(a)=a\circ \bc^{(k)}_R$, respectively.
If $b:R'\ra S$ is another local morphisms such that $c^{(k)}(a)=d^{(k)}(b)$ then the unique local map $a;_kb: R\otimes_kR'\ra S$ such, that $a;_kb\circ \kappa^1= a$ and $a;_kb\circ \kappa^2= b$, is the composition of $a$ and $b$ in $S^*$. For more details on functor $(-)^*$ see \cite{Z2}.

Now we will set up the adjunction
\begin{center} \xext=2500 \yext=300
\begin{picture}(\xext,\yext)(\xoff,\yoff)
\putmorphism(0,200)(1,0)[\phantom{Mod_\otimes((\ofs)^{op},Set)}`
\phantom{\mComp}`\widetilde{(-)}]{2500}{1}a
\putmorphism(0,100)(1,0)[\phantom{Mod_\otimes(\ofs^{op},Set)}`
\phantom{\mComp}`\widehat{(-)}=\mComp((\simeq)^*, - )]{2500}{-1}b
\putmorphism(0,150)(1,0)[Mod_\otimes(\ofs^{op},Set)`\mComp`]{2500}{0}b
\end{picture}
\end{center}
which will turn out to be an equivalence of categories. The
functor $\widehat{(-)}$ is sending a many-to-one computad $\cP$ to
a functor
\[ \widehat{\cP} = \mComp((-)^*,\cP) : \ofs^{op} \lra Set \]
$\widehat{(-)}$ is defined on morphism in the obvious way, by
composition. We have

\begin{lemma} Let $\cP$ be a many-to-one computad. Then $\widehat{\cP}$
defined above sends tensor squares to pullbacks.
\end{lemma}

{\it Proof.}~This is an immediate consequence of the fact that the
functor \mbox{$(-)^*:\ofs\ra\mComp$} sends tensor squares to pushouts, Corollary 13.3 in
[Z2]. $~~\Box$

The functor $\widetilde{(-)}$ that we describe below is the induced functor described
in Proposition \ref{GTmodel} for the model $(-)^*:\ofs^{op}\ra (\mComp)^{op}$.
As we need to establish some properties of $\widetilde{(-)}$ we give here a more detailed  description.

Suppose we have a model $F : \ofs^{op} \lra Set$. We shall define a many-to-one computad
$\widetilde{F}$. As the set of $n$-cells of $\widetilde{F}$ we take
\[ \widetilde{F}_n = \coprod_S F(S) \]
where the coproduct\footnote{ In fact, we think about such a
coproduct $\coprod_S F(S)$ as if it were to be taken over
sufficiently large (so that each isomorphism type of ordered face
structures is represented) set of ordered face structures $S$ of
dimension at most $n$. Then, if ordered face structures $S$ and
$S'$ are isomorphic via (necessarily unique) monotone isomorphism
$h:S'\ra S$, then the cells $x\in F(S)$ and $x'\in F(S')$ are
considered equal iff $F(h)(x)=x'$.} is taken over all (up to a
monotone isomorphisms) ordered face structures $S$ of dimension at
most $n$. By $\kappa^{F,S}_n:F(S)\lra\widetilde{F}_n$ we denote
the coprojection into the coproduct.
For $k\leq n$, the identity map
\[ 1^{(n)} : \widetilde{F}_k \lra \widetilde{F}_n \]
is the obvious embedding induced by identity maps on the
components of the coproducts. For $k\leq n$, we define the
$k$-domain and the $k$-codomain functions in $\widetilde{F}$
\[ d^{(k)},c^{(k)} : \widetilde{F}_n \lra \widetilde{F}_k. \]
Abstractly, $d^{(k)}$ is the unique map, that makes
the diagram
\begin{center} \xext=1000 \yext=650
\begin{picture}(\xext,\yext)(\xoff,\yoff)
 \setsqparms[1`-1`-1`1;1000`500]
 \putsquare(0,100)[\widetilde{F}_n`\widetilde{F}_k`F(S)`F(\bd^{(k)}S);d^{(k)}`\kappa^{F,S}_n`\kappa^{F,d^{(k)}S}_n`F(\bd^{(k)}_S)]
\end{picture}
\end{center}
commute, for any ordered face structure $S$. $c^{(k)}$ is defined
similarly. In more concrete terms  $d^{(k)}$ and $c^{(k)}$ are
defined as follows. Let $S$ be an ordered face structure of
dimension at most $n$, $a\in F(S)\lra \widetilde{F}_n$ an $n$-cell
in $\widetilde{F}$. We have in $\ofs$ the morphisms of the $k$-th
domain and the $k$-th codomain introduced in \cite{Z2}:
\begin{center}
\xext=800 \yext=400 \adjust[`I;I`;I`;`I]
\begin{picture}(\xext,\yext)(\xoff,\yoff)
 \settriparms[-1`-1`0;400]
 \putAtriangle(0,0)[S`\bd^{(k)} S`\bc^{(k)} S;\bd^{(k)}_S`\bc^{(k)}_S`]
\end{picture}
\end{center}
We put
\[ d^{(k)}(a) = F(\bd^{(k)}_S)(a)\in F(\bd^{(k)}S) \lra\widetilde{F}_k ,\]
\[ c^{(k)}(a) = F(\bc^{(k)}_S)(a)\in F(\bc^{(k)}S) \lra\widetilde{F}_k. \]

Finally, we define the compositions in $\widetilde{F}$. Again we
shall do it first abstractly and then in concrete terms.  Note
that the pullback
\begin{center} \xext=800 \yext=550
\begin{picture}(\xext,\yext)(\xoff,\yoff)
 \setsqparms[1`1`1`1;800`400]
 \putsquare(0,50)[\widetilde{F}_n\times_{\widetilde{F}_k}\widetilde{F}_n`\widetilde{F}_n`\widetilde{F}_n`\widetilde{F}_k;
 \pi_1`\pi_0`d^{(k)}`c^{(k)}]
\end{picture}
\end{center}
can be describe as a coproduct
\[ \widetilde{F}_n\times_{\widetilde{F}_k}\widetilde{F}_n
=\coprod_{S,S'} F(S)\times_{F(\bc^{(k)}S)}F(S')\hskip 5mm
(\cong\coprod_{S,S'} F(S\otimes_kS')) \] where the coproduct is
taken over all (up to monotone isomorphisms) pairs of ordered face
structures $S$ and $S'$ of dimension at most $n$ such that
$\bc^{(k)}S=\bd^{(k)}S'$. The coprojections are denoted by
\[ \kappa^{F,S,S'}_{n,k}:F(S)\times_{F(\bc^{(k)}S)}F(S')\lra \coprod_{S,S'}
F(S)\times_{F(\bc^{(k)}S)}F(S')=\widetilde{F}_n\times_{\widetilde{F}_k}\widetilde{F}_n\]
Then the composition morphism
\[ \textbf{;}_k : \widetilde{F}_n\times_{\widetilde{F}_k}\widetilde{F}_n\lra
\widetilde{F}_n\]
 is the unique map that for any pair $S$, $S'$ as above
makes the square
\begin{center} \xext=1200 \yext=750
\begin{picture}(\xext,\yext)(\xoff,\yoff)
 \setsqparms[1`-1`-1`1;1200`600]
 \putsquare(0,50)[\widetilde{F}_n\times_{\widetilde{F}_k}\widetilde{F}_n`\widetilde{F}_n`F(S)\times_{F(\bd^{(k)}S)}F(S')`
 F(S\otimes_kS');
 \textbf{;}_k`\kappa^{F,S,S'}_{n,k}`\kappa^{F,S\otimes_kS'}_{n}`\zeta_{S,S'}]
\end{picture}
\end{center}
commute, where $\zeta_{S,S'}$ is the inverse of the canonical
isomorphism
\[F(S\otimes_kS')\lra F(S)\times_{F(\bc^{(k)}S)}F(S')\]
that exists as $F$ preserves special pullbacks. In concrete terms, the composition in
$\widetilde{F}$ can be described as follows.  Let $k<n$,
$dim(S),dim(S')\leq n$, $\bc^{(k)}S=\bd^{(k)}S'$,
\[ a\in F(S) \lra\widetilde{F}_n \hskip 10mm b\in F(S')
\lra\widetilde{F}_n,\]
 such that
 \[ c^{(k)}(a)=F(\bc^{(k)}_S)(a)=F(\bd^{(k)}_{S'})(b) = d^{(k)}(b). \]
We shall define the cell $a\textbf{;}_kb\in\widetilde{F}_n$. We
have a tensor square in $\ofs$:
\begin{center} \xext=1500 \yext=680
\begin{picture}(\xext,\yext)(\xoff,\yoff)
 \setsqparms[1`-1`-1`1;800`500]
 \putsquare(0,100)[S`S\otimes_kS'`{\bc^{(k)}S}`S';
 \kappa_S`\bc^{(k)}_S`\kappa_{S'}`\bd^{(k)}_{S'}]
\end{picture}
\end{center}
As $F$ is a model of $\ofs$  the square
\begin{center} \xext=1500 \yext=670
\begin{picture}(\xext,\yext)(\xoff,\yoff)
 \setsqparms[1`-1`-1`1;800`500]
 \putsquare(0,100)[F(S)`F(S\otimes_kS')`F({\bc^{(k)}S})`F(S');
 F(\kappa_S)`F(\bc^{(k)}_S)`F(\kappa_{S'})`F(\bd^{(k)}_{S'})]
\end{picture}
\end{center}
is a pullback in $Set$. Thus there is a unique element
\[ x\in F(S\otimes_kS')\lra \widetilde{F}_n \]
such that
\[ F(\kappa_S)(x)=a, \hskip 10mm   F(\kappa_{S'})(x)=b. \]
 We put
\[ a\textbf{;}_kb = x. \]
This ends the definition of $\widetilde{F}$.

 For a morphism $\alpha :F\lra G$ in $Mod_\otimes(\ofs^{op},Set)$ we put
 \[ \widetilde{\alpha} = \{ \widetilde{\alpha}_n :
 \widetilde{F}_n\lra \widetilde{G}_n \}_{n\in\o} \]
 such that
 \[ \widetilde{\alpha}_n = \coprod_S \alpha_S : \widetilde{F}_n\lra \widetilde{G}_n\]
where the coproduct is taken over all (up to monotone isomorphism)
ordered face structures $S$ of dimension at most $n$. This ends
the definition of the functor $\widetilde{(-)}$.

We have

\begin{proposition}
The functor
\[ \widetilde{(-)} : Mod_\otimes(\ofs^{op},Set)\lra \mComp \]
is well defined.
\end{proposition}

{\it Proof.}~The verification that $\widetilde{(-)}$ is a functor
into $\oC$ is left for the reader.  We shall verify that, for model  $F:{\ofs}^{op}\lra
Set$ of $\ofs$, $\widetilde{F}$ is a many-to-one computad, whose $n$-indets
are
 \[ |\widetilde{F}|_n = \coprod_{P\in \pfs, dim(P)=n} F(P)\; \lra \;
 \coprod_{S\in \ofs, dim(S)\leq n} F(S) = \widetilde{F}_n.\]

Let $\cP$ be the $n$-truncation of $\widetilde{F}$ in  $\mnComma$,
i.e. $\cP=\widetilde{F}^{\natural,n}$ in the notation from
Appendix of [Z2]. We shall show that $\widetilde{F}_n$ is in a
bijective correspondence with $\overline{\cP}_n$ described in
[Z2]. We define a function
\[ \varphi : \overline{\cP}_n \lra \widetilde{F}_n\]
so that for a cell $f:S^{\sharp,n}\lra \cP$ in $\overline{\cP}_n$
we put
\[  \varphi(f) = \left\{ \begin{array}{ll}
        1_{f_{n-1}(S)}   & \mbox{if $dim(S)<n$,}  \\
        f_{n}(m_S)   & \mbox{if $dim(S)=n$, $S$ principal, $S_n=\{ m_S \}$}  \\
        \varphi(f^{\da \breve{a}});_k \varphi(f^{\ua \check{a}})& \mbox{if $dim(S)=n$, $\breve{a}\in Sd(S)_k$.}
                                    \end{array}
                \right. \]
and the morphisms in $\varphi(f^{\da \check{a}})$ and
$\varphi(f^{\ua \check{a}})$ in $\mnComma$ are obtained by
compositions so that the diagram
\begin{center}
\xext=1000 \yext=700
\begin{picture}(\xext,\yext)(\xoff,\yoff)
\settriparms[0`1`-1;300]
 \putDtriangle(0,0)[(S^{\ua \check{a}})^{\sharp,n}`S^{\sharp,n}`(S^{\da \check{a}})^{\sharp,n};``]
 \putmorphism(400,300)(1,0)[`\cP`f]{600}{1}a
 \putmorphism(0,700)(3,-1)[``f^{\ua \check{a}}]{1100}{1}r
 \putmorphism(400,50)(3,1)[``f^{\da \check{a}}]{280}{1}r
\end{picture}
\end{center}
commutes. We need to verify, by induction on $n$, that $ \varphi$
is well defined, bijective and that it preserves compositions,
identities, domains, and codomains.

We shall only verify (partially) that $\phi$ is well defined, i.e.
that the definition of $\phi$ for any  non-principal ordered face
structure $S$ of dimension $n$ does not depend on the choice of
the saddle point of $S$. Let $\check{a},\check{b}\in Sd(S)$. We
shall show, in case $dim(a)=dim(b)=k$ and $a<_l b$, that we have
 \[ \varphi(f^{\da \check{a}});_k \varphi(f^{\ua \check{a}})=
 \varphi(f^{\da \check{b}});_k \varphi(f^{\ua \check{b}})\]

Using  Lemma 12.6 of [Z2] 
and the fact that $(-)^{\sharp,n}$ preserves special pushouts
(Corollary 13.2 of [Z2]), we have
\begin{eqnarray*}
 \varphi(f^{\da \check{a}});_k \varphi(f^{\ua \check{a}}) = \\
 = \varphi(f^{\da a});_k (\varphi(f^{\ua \check{a}\da b});_k \varphi(f^{\ua \check{a}\ua \check{b}})) = \\
 = (\varphi(f^{\da a});_k \varphi(f^{\ua \check{a}\da b}));_k \varphi(f^{\ua \check{b}\ua \check{a}}) =\\
  = \varphi(f^{\da a};_k f^{\ua \check{a}\da b});_k \varphi(f^{\ua \check{b}\ua \check{a}}) = \\
  = \varphi(f^{\da b};_k f^{\ua \check{b}\da a});_k \varphi(f^{\ua \check{b}\ua \check{a}}) =  \\
  = (\varphi(f^{\da b});_k \varphi(f^{\ua \check{b}\da a}));_k \varphi(f^{\ua \check{b}\ua \check{a}}) = \\
  = \varphi(f^{\da b});_k (\varphi(f^{\ua \check{b}\da a});_k \varphi(f^{\ua \check{b}\ua \check{a}})) =  \\
 = \varphi(f^{\da b});_k \varphi(f^{\ua \check{b}})
\end{eqnarray*}
The reader can compare these calculations with the those, in the
same case, of Proposition 13.1 of [Z2] 
($f$ replaces $\varphi$ and $\varphi$ replaces $F$). So there is
no point to repeat the other calculations here. $~~\Box$

 For $\cP$ in $\mComp$ we define a computad map
\[ \eta_\cP : \cP\lra \widetilde{\widehat{\cP}} \]
so that for $x\in \cP_n$ we put
\begin{center} \xext=1600 \yext=150
\begin{picture}(\xext,\yext)(\xoff,\yoff)
 \putmorphism(0,30)(1,0)[\eta_{\cP,n}(x)=\tau_x : T_x^*`\cP\in \widehat{\cP}(T_x)`]{800}{1}b
 \putmorphism(800,30)(1,0)[\phantom{\cP\in
 \widehat{\cP}(T_x)}`\widetilde{\widehat{\cP}}_n`\kappa_n^{T_x}]{800}{1}a
\end{picture}
\end{center}
such that $\tau_x(1_{T_x})=x$.

For $F$ in $Mod_\otimes(\ofs^{op},Set)$ we define a natural
transformation
\[ \varepsilon_F : \widehat{\widetilde{F}} \lra F, \]
such that, for an ordered face structure $S$ of dimension $n$,
\[ (\varepsilon_F)_S : \widehat{\widetilde{F}}(S) \lra F(S) \]
and $g : S^*\ra \widetilde{F} \in \widehat{\widetilde{F}}(S)$ we
put
\[ (\varepsilon_F)_S(g) = g_n(1_S).\]

\begin{proposition}\label{presheaf_eq3}
The functors
\begin{center} \xext=2500 \yext=300
\begin{picture}(\xext,\yext)(\xoff,\yoff)
\putmorphism(0,200)(1,0)[\phantom{Mod_\otimes(\ofs^{op},Set)}`
\phantom{\mComp}`\widetilde{(-)}]{2500}{1}a
\putmorphism(0,100)(1,0)[\phantom{Mod_\otimes(\ofs^{op},Set)}`
\phantom{\mComp}`\widehat{(-)}=\mComp((\simeq)^*, - )]{2500}{-1}b
\putmorphism(0,150)(1,0)[Mod_\otimes(\ofs^{op},Set)`\mComp`]{2500}{0}b
\end{picture}
\end{center}
together with the natural transformations $\eta$ and $\varepsilon$
defined above form an adjunction
($\widehat{(-)}\dashv\widetilde{(-)}$). It establishes an
equivalence of categories $Mod_\otimes(\ofs^{op},Set)$ and $\mComp$.
\end{proposition}

{\it Proof.}~ The fact that both $\eta$ and $\varepsilon$ are
bijective on each component follows immediately from Proposition
15.1 of [Z2]. So we shall verify the triangular equalities only.

Let $\cP$ be a computad, and $F$ be a functor in
$Mod_\otimes(\ofs^{op},Set)$. We need to show that the triangles
\begin{center}
\xext=2000 \yext=400 \adjust[`I;I`;I`;`I]
\begin{picture}(\xext,\yext)(\xoff,\yoff)
 \settriparms[-1`1`1;400]
 \putAtriangle(0,0)[\widehat{\widetilde{\widehat{\cP}}}`\widehat{\cP}`\widehat{\cP};
 \widehat{\eta_\cP}`\varepsilon_{\widehat{\cP}}`1_{\widehat{\cP}}]
 \putAtriangle(1200,0)[\widetilde{\widehat{\widetilde{F}}}`\widetilde{F}`\widetilde{F};
 \eta_{\widetilde{F}}`\widetilde{\varepsilon_{F}}`1_{\widetilde{F}}]
\end{picture}
\end{center}
commute. So let $f:S^*\ra \cP \in \widehat{\cP}(S)$. Then, we have
\[ \varepsilon_{\widehat{\cP}}\circ \widehat{\eta_\cP} (f) =
\varepsilon_{\widehat{\cP}}(\eta_\cP\circ f)= (\eta_\cP\circ
f)_n(1_S) = \]
\[ =   (\eta_\cP)_n(f_n(1_S)) = \tau_{f_n(1_S)} = f \]
Last equation follows from the fact that $(\tau_{f_n(1_S)})_n(1_S)
= f_n(1_S)$ and Proposition 15.1 of [Z2]. Now let $x \in F(S)\lra
\widetilde{F}_n$. Then we have

\[ \widetilde{\varepsilon_F}\circ \eta_{\widetilde{F}}(x)=
\widetilde{\varepsilon_F}(\tau_x) = (\tau_x)_n(1_{T_x}) =x  \] So
both triangles commute, as required. $~~\Box$

If we compose the three established adjoint equivalences we get from
Propositions \ref{presheaf_eq1},  \ref{presheaf_eq2}, and \ref{presheaf_eq3}

\begin{corollary}
The functor
\[ \widehat{(-)} :  \mComp \lra Set^{\pfs^{op}} \]
such that for a many-to-one computad $X$,
\[ \widehat{X} = \mComp((-)^*,X): \pfs^{op}\lra Set  \]
is an equivalence of categories.
\end{corollary}

The fact that the category $\mComp$ is a presheaf category  was
first established in \cite{HMZ} using an earlier result from
\cite{HMP}.  From this we know that the category of $\mComp$ is
equivalent the category of presheaves on the category of
multitopes $\bf Mlt$ introduced in \cite{HMP}.

\begin{theorem}
The category $\pfs$ of principal ordered face structure is
equivalent to the category of multitopes  $\bf Mlt$.
\end{theorem}

{\it Proof.}~ The categories of presheaves on both categories are
equivalent to the category of many-to-one computads. As these
categories have no nontrivial idempotents they must be equivalent.
 $~\Box$

\section{The shapes of cells in computads}\label{shape}

Let $\Comp^{?/?}$ be a full subcategory of the category of computads $\Comp$\footnote{By a computad we mean here an $\o$-category $C$ that is
levelwise free, i.e. if we truncate it to an $n+1$-category $C_{n+1}$ then it arises as an $n$-category $C_n$ with freely added $n+1$-indeterminate
cells (=indets). Morphisms of computads are $\o$-functors that are required to send indets to indets.} of some kind of computads.
The particular examples we have in mind and we will be refereing to later are free categories over graphs $\Comp^{1/1}_1$, one-to-one
computads $\Comp^{1/1}$, positive-to-one computads $\pComp$, many-to-one computads $\mComp$, or even all computads $\Comp$.

One of the ways to think about the shape of a cell $\alpha$ in a computad $C$ from the category $\Comp^{?/?}$ is the following. We
consider the category of pointed computads $\Comp^{?/?}_*$ whose objects are computads with chosen cells and morphisms are
computad maps preserving the distinguished cells. Then the {\em shape} of a cell $\alpha\in C$ (if exists) can be identified with
the initial object of the slice category $\Comp_*\da (C,\alpha)$. It is obvious that the computad maps preserve so understood
shapes i.e. if $f:C\ra D$ is a computad map, $\tau_\alpha:(S,s)\ra (C,\alpha)$ is the initial object of $\Comp_*\da (C,\alpha)$ then
$F\circ \tau_\alpha:(S,m)\ra (D,f(\alpha))$ is the initial object of $\Comp_*\da (D,f(\alpha))$. Unfortunately not every cell has a
shape.  For example if we take two 2-indets $\alpha$ and $\beta$ whose domain and codomain is $1_x$ the identity of a $0$-cell $x$
then $\beta\circ_0\alpha$ does not have a shape. This 'innocent' problem is responsible for very serious complications and it is
one of the reasons for the restriction of shapes of cells in weak $\o$-categories to more manageable shapes like one-to-one,
many-to-one, etc. Note that the shape $(S,m)$ of the cell $\alpha$ is not necessarily determined by what we can call the (pure) shape
$S$.

Now assume that all cells in all computads in the given category of computads $\Comp^{?/?}$ have shapes.
To define the category $\Shape^{?/?}$ of shapes of cells for $\Comp^{?/?}$  we could just take all the computad maps between shapes of all cells in computads from  $\Comp^{?/?}$. But such morphisms can identify different shapes by making them isomorphic.  This is why we shall take a longer route by specifying some of the morphisms that we definitely want in the category $\Shape^{?/?}$ and then we shall generate all the other morphisms inside $\Comp^{?/?}$
via composition and graded tensor operation. The closure under the later operation is to ensure that the graded tensor operation is functorial.

First kind of morphisms we shall consider comes from the fact that we have in computads the $k$-domain $d^{(k)}$ and the $k$-codomain $c^{(k)}$ operations that associate the domain and the codomain of dimension $k$, respectively. Thus if $(S,m)$ is a shape, $m$ is a cell in $S$ of dimension $n$ and $k\leq n$, then $(S,d^{(k)}(m))$ and $(S,c^{(k)}(m))$ are  pointed computads in $Comp^{?/?}$. Thus the cells $d^{(k)}(m)$ and  $c^{(k)}(m)$ in $S$ have shapes which we denote
\[ \bd^{(k)}_S: (\bd^{(k)}S,d^{(k)}(m))\lra  (S,d^{(k)}(m)) \]
\[ \bc^{(k)}_S: (\bc^{(k)}S,c^{(k)}(m))\lra  (S,c^{(k)}(m)) \]
In particular, we have the computad maps
 \[ \bd^{(k)}_S: \bd^{(k)}S\lra S,\hskip5mm  \bc^{(k)}_S: \bc^{(k)}S\lra S\]
such that $\bd^{(k)}_S(d^{(k)}(m))=d^{(k)}(m)$ and $\bc^{(k)}_S(c^{(k)}(m))=c^{(k)}(m)$. Note that both $d^{(k)}(m)$ and $c^{(k)}(m)$ name two different cells in two different computads.

The second kind of morphism comes from the fact that we can
(de)compose cells in computads. Suppose $(S,m)$ is a shape such that the cell $m$ can be decomposed as $m=m_1;_km_2$.
Then we have shapes of $m_1$ and $m_2$ in $S$:
\[ \kappa^1:(S_1,m_1)\lra (S,m_1),\hskip 5mm \kappa^2:(S_1,m_2)\lra (S,m_2) \]
so that $m=m_1;_km_2=\kappa^1(m_1);_k\kappa^2(m_2)$ in $S$.
Here again both $m_1$ and $m_2$ name two different cells in two different computads.
If we denote by $(S_3,m_3)$ the shape of $c^{(k)}(m_1)=d^{(k)}(m_2)$ in $S$ then
we obtain a commuting square
 \begin{center}
\xext=800 \yext=600
\begin{picture}(\xext,\yext)(\xoff,\yoff)
 \setsqparms[1`-1`-1`1;800`500]
 \putsquare(0,50)[(S_1,m_1)`(S,m)`(S_3,m_3)`(S_2,m_2);\kappa^1`\bd^{(k)}_{S_1}`\kappa^2`\bc^{(k)}_{S_2}]
 \end{picture}
\end{center}
called the tensor square.

Note that it is very likely that the computad $S$
will turn out to be the pushout $S_1+_{S_3}S_2$ in $Comp^{?/?}$ but this doesn't mean that $(S,m)$ will be the pushout  $(S_1,m_1)+_{(S_3,m_3}(S_2,m_2)$ in  $Shape^{?/?}$. 

The graded tensor operation is defined as follows. Suppose we have tensor squares defined form decompositions of cells $m=m_1;_km_2$ and $m'=m'_1;_km'_2$
in shapes $(S,m)$ and $(S',m')$, respectively, and for some morphisms $f_1$, $f_2$, $f_3$,  the squares
\begin{center}
\xext=2000 \yext=700
\begin{picture}(\xext,\yext)(\xoff,\yoff)
 \setsqparms[-1`1`1`-1;1000`500]
 \putsquare(0,100)[(S_1,m_1)`(S_3,m_3 )`(S'_1,m'_1)`(S'_3,m'_3);\bc^{(k)}_{S_1}`f_1`f_3`\bc^{(k)}_{S'_1}]
 \setsqparms[1`0`1`1;1000`500]
 \putsquare(1000,100)[\phantom{(S_3,m_3)}`(S_2,m_2)`\phantom{(S'_3,m'_3 )}`(S'_2,m'_2);\bd^{(k)}_{S_2}``f_2`\bd^{(k)}_{S'_2}]
 \end{picture}
\end{center}
commute. Then we require to exist a unique morphism $f_1\otimes_kf_1: S\ra S'$, called a {\em graded tensor} of $f_1$ and $f_1$, making the squares
\begin{center}
\xext=2000 \yext=700
\begin{picture}(\xext,\yext)(\xoff,\yoff)
 \setsqparms[1`1`1`1;1000`500]
 \putsquare(0,100)[(S_1,m_1)`(S,m)`(S'_1,m'_1)`(S',m' );\kappa^1_{S_1}`f_1``\kappa^1_{S'_1}]
 \setsqparms[-1`0`1`-1;1000`500]
 \putsquare(1000,100)[\phantom{(S,m)}`(S_2,m_2)`\phantom{(S',m')}`(S'_2,m'_2);\kappa^2_{S_2}`f_1\otimes_kf_2`f_2`\kappa^2_{S'_2}]
 \end{picture}
\end{center}
commute.

The category of $?/?$-shapes $Shape^{?/?}$ has as objects shapes $(S,m)$ of cells in $Comp^{?/?}$ and as morphisms
the least class of computad morphism containing $\bd^{(k)}$, $\bd^{(k)}$, $\kappa^1$, $\kappa^2$, identities closed under composition and graded tensor operation.

\vskip 2mm
{\em Examples.}

1. For the category of free categories over graphs i.e. $Comp^{1/1}_1$ the category of shapes  defined above in (equivalent to) $\Delta_0$. Recall that $\Delta_0$ is the full subcategory of the category of graphs whose objects are  linear graphs $[n]$ with $n$ edges and $n+1$ vertices. The morphism $\bd^{(0)}_{[n]}:[0]\ra [n]$ is the inclusion sending the unique vertex of $[0]$ to the first vertex of $[n]$ and  $\bc^{(0)}_{[n]}:[0]\ra [n]$ is the inclusion sending the unique vertex of $[0]$ to the last vertex of $[n]$. The morphisms $\kappa^1_{[n]}:[n]\ra [n+m]$ and $\kappa^2_{[m]}:[n]\ra [n+m]$ are the inclusions onto the first $n$ and the last $m$ edges of the graph $[n+m]$, respectively. The tensor morphisms are also obvious. Note that in this case we generate all the graph morphisms between the objects of $\Delta_0$.

2. The shapes of cells for the category of one-to-one computads $Comp^{1/1}$ are determined by what is called in different terminologies globular cardinals, simple $\o$-graphs $s\oG$,  $T$-cardinals for the free category monad on $\o$-graphs.  By this I mean that for every $n$-cell $\alpha$ in every one-to-one computad $C$ there is a unique (up to iso) simple $\o$-graph $S$ and a unique cell $S$ in the $\o$-category $S^*$ generated\footnote{The $n$-cells of $S^*$ can be identified with simple subgraphs of $S$ of dimension at most $n$} by $S$, and a unique pointed computad morphism $\tau_\alpha:(S^*,S)\ra (C,\alpha)$. Moreover this map $\tau_\alpha$ is the initial object in $Comp^{1/1}_*\da (C,\alpha)$. The shape $(S^*,S)$ is uniquely determined by $\o$-graph $S$ and even by the $\o$-category $S^*$.

3. The shapes of cells for the category of positive-to-one computads $\pComp$, see \cite{Z1}, are determined by positive face structures. The category $\fs$
of positive face structures is the category of shapes for $\pComp$. Despite the fact that it is considerably more complicated than $s\oG$, it shares some good properties of $s\oG$. For example the embedding $(-)^*:\fs \ra \pComp$ is full.

4. The shapes of cells for the category of many-to-one computads $\mComp$, see \cite{Z2}, are determined by ordered face structures. The category $\ofs$ of positive face structures and monotone maps is the category of shapes for $\mComp$. Here however the theory changes considerably. The main reason is that the $(-)^*:\ofs \ra \mComp$ is not full. The full image of $(-)^*$ in this case is the category of ordered face structures and local maps.

5. As not all the cells in arbitrary computads have shape, there is no category of shapes for the category of all computads $Comp$.

\section{Pra monads and nerves}\label{pra monads}


The idea that algebras can be presented as a full subcategory of the category of free algebras preserving
some limits goes back to the thesis of our jubilee. In \cite{Lawvere} F.W. Lawvere  have shown that finitary algebras
can be presented as a full subcategory of presheaves on the finitely generated free algebras that preserve
some finite products. The next step was made by F.E.Linton, c.f. \cite{Linton}, when he has shown that
for any\footnote{In fact F.E.Linton needed some minor size restricting condition called {\em tractable} saying that
operation of any (possibly infinitary) arity form a set.} monad
$T$ on $Set$ it is true that the category $Alg(T)$ of
the Eilenberg-Moore algebras for $T$ is equivalent to the category of product preserving functors from the dual of the category $\cK(T)$ of Kleisli
algebras\footnote{At that time it was not expressed in these terms.}. He also noticed that under further size restrictions
on $T$ one can take an essentially small full subcategory of $\cK(T)$.

The recent development due to T. Leinster, c.f. \cite{Leinster1}, and then to M. Weber, c.f. \cite{Weber}, brought some new light on this construction.
Below I will describe briefly the theory developed by them but not in full generality of M. Weber and changing slightly the perspective occasionally.

T. Leinster's setup consists of a parametric right adjoint monad $(T,\nu,\mu)$, pra monad for short, on a presheaf category $Set^{\cC^{op}}$.
By this he means that both natural transformations are cartesian and that functor $T$  a parametric right adjoint i.e. that the functor \mbox{$T_1:Set^{\cC^{op}}\lra Set^{\cC^{op}}\da T(1)$} induced by $T$ has a left adjoint.
He has shown that in this case  there is a canonical choice for a small category $\theta_T$ a full subcategory of $\cK(T)$ and a canonical
choice of the limits in  $\theta_T$ so that the category of presheaves preserving those limits is equivalent to $Alg(T)$.

A functor $T$ defined on a presheaf category is pra iff it preserves wide pullbacks iff it is family representable, c.f. \cite{Leinster1},  \cite{Weber}.
Recall that a functor on a presheaf category $Set^{\cC^{op}}$ is a family representable iff for every object $c\in \cC$ there is a set of objects $\{ T_{c,i}\}_{i\in I_c}$ of $Set^{\cC^{op}}$ such that we have an isomorphism of functors
$ev_c \circ T \cong \coprod_{i\in I_c}  Y(T_{c,i})$, where $ev_c:Set^{\cC^{op}}\ra Set$ is the evaluation on $c$, and $Y(T_{c,i})$ is the covariant functor  representable by $T_{c,i}$. The category $\theta_T$ is the full subcategory of $Alg(T)$ whose objects are the free $T$-algebras over the representing objects
$\{ T_{c,i}\}_{i\in I_c, c\in \cC}$.

In M. Weber's terminology the objects of form $T_{c,i}$ for $c\in C$, $i\in I_c$ are called $T${\em -cardinals}. In order to make the distinction I will call his category $\Theta_0$ a full subcategory of $Set^{\cC^{op}}$  as {\em the category of $T$-cardinals} and the full image of it in $Alg(T)$ denoted by him $\Theta_T$ as the {\em $T$-cardinal algebras}.

M. Weber is considering a more general setup than T. Leinster. The monad $(T,\nu,\mu)$ is defined on a cocomplete category $\cA$. In this more general situation the choice of the category of $T$-cardinals, called there the category of arities, does not need to be canonical and is given explicitly as full dense subcategory of $\cA$.  M. Weber also requires $\eta$ and $\mu$ to be cartesian but the condition on $T$ is slightly more technical and I will not recall it here. The more general setup covers some cases not
covered by T. Leinster approach but the additional level of generality  has in the present context only restricted and negative application to which I will come back later. On the other hand, the Theorem 4.10 in \cite{Weber} seems to be more informative even when applied to the original setup of T. Leinster.
In fact I will state it in combination with other results from \cite{Weber} and \cite{Leinster1} in the form that is relevant to the present context.

Now let $p\Theta$ be a small category  $(T,\eta,\mu)$ be a pra monad on a presheaf category $Set^{p\Theta}$, $\Theta_0$ the full subcategory of $Set^{p\Theta}$ whose objects are $T$-cardinals, $\Theta_T$ the full image of the category $\Theta_0$ in $Alg(T)$,  $\cM$ the class of morphisms in $\Theta_0$. Then, following M. Weber, we can conclude that there is a class $\cE$ orthogonal to $\cM$ so that $(\cE,\cM)$ form a factorization system in $\Theta_T$ moreover we have a commuting square of categories and functors
 \begin{center}
\xext=800 \yext=600
\begin{picture}(\xext,\yext)(\xoff,\yoff)
 \setsqparms[1`1`1`1;1000`500]
 \putsquare(0,50)[Alg(T)`Set^{\Theta^{op}_T}`Set^{p\Theta^{op}}`Set^{\Theta^{op}_0};\cN_T`U`i^*`]
 \end{picture}
\end{center}
which is a pseudo-pullback, where the horizontal maps are the obvious maps generated by the (full) embeddings   $\Theta_0\lra Set^{p\Theta}$ and
$\Theta_T\lra Alg(T)$. As the first embedding is full (and faithful) so is the nerve functor $\cN_T$. We can think about this result as saying that
all the monadic functors (like $U:Alg(T)\lra Set^{p\Theta^{op}}$)  for pra monads can be obtained via pseudo-pullbacks along full and faithful functors from particularly simple monadic functors namely those coming from presheaf pra monads, (like $i^*: Set^{\Theta^{op}_T}\lra Set^{\Theta^{op}_0}$), see Example 1 below.

\vskip 2mm
{\em Examples.}

1. There is a whole class of simple examples of pra monads. Let $\Xi$ be a small category with a factorization system $(\cE,\cM)$\
such that for any morphisms $e,e_1,e_2,e'\in \cE$ and $m,m_2\in \cM$ making the left hand diagram below commutes there exist $e'_1,e'_2\in \cE$ and $m_1\in\cM$, unique up to an isomorphism,  making the right hand diagram below commutes
 \begin{center}
\xext=2400 \yext=780
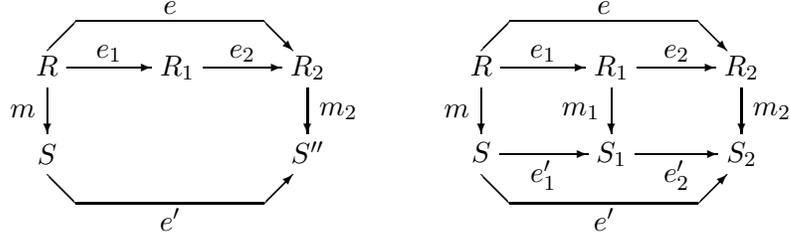
\begin{picture}(\xext,\yext)(\xoff,\yoff)
 \setsqparms[1`1`0`0;450`300]
 \putsquare(0,250)[R`R_1`S`;e_1`m``]
 \setsqparms[1`0`1`0;450`300]
 \putsquare(450,250)[\phantom{R_1}`R_2`\phantom{S_1}`S'';e_2``m_2`]

 \put(750,80){\vector(1,1){100}}
 \put(100,80){\line(-1,1){100}}
\put(100,80){\line(1,0){650}}

\put(400,0){\makebox(50,50){$e'$}}
 \put(750,710){\vector(1,-1){100}}
 \put(100,710){\line(-1,-1){100}}
\put(100,710){\line(1,0){650}}

\put(400,730){\makebox(50,50){$e$}}


 \setsqparms[1`1`1`1;450`300]
 \putsquare(1500,250)[R`R_1`S`S_1;e_1`m``e_1']
 \setsqparms[1`0`1`1;450`300]
 \putsquare(1950,250)[\phantom{R_1}`R_2`\phantom{S_1}`S_2;e_2`m_1`m_2`e_2']

 \put(2250,80){\vector(1,1){100}}
 \put(1600,80){\line(-1,1){100}}
\put(1600,80){\line(1,0){650}}

\put(1900,0){\makebox(50,50){$e'$}}
 \put(2250,710){\vector(1,-1){100}}
 \put(1600,710){\line(-1,-1){100}}
\put(1600,710){\line(1,0){650}}

\put(1900,730){\makebox(50,50){$e$}}
 \end{picture}
\end{center}
Let $\Xi_\cM$ be the subcategory for $\Xi$ consisting of all objects of $\Xi$ and morphisms from the class $\cM$ only. We have
a non-full surjective on objects embedding $i:\Xi_\cM\ra \Xi$. Then the functor $i^*:Set^{\Xi^{op}}\ra Set^{\Xi^{op}_\cM}$ is a monadic and the monad ($(T^i,\eta^i,\mu^i)$ on $Set^{\Xi^{op}_\cM}$ induced by $i$ is pra. For $X\in Set^{\Xi^{op}_\cM}$ and $S\in \Xi_\cM$ the functor $T^i$ is given by
\[ T^i(X)(S)=\coprod_{e:S\ra R} X(R)\]
where the coproduct is taken over all (up to iso) morphisms in $\cE$ with domain $S$. The additional condition that we require for the factorization system
is needed to show that the natural transformation $\mu^i$ is cartesian. As for such monad $T^i$ not only the base category $Set^{\Xi^{op}_\cM}$ but also the category of algebras $Alg(T)$ is a presheaf category $Set^{\Xi^{op}}$, I will call such monads {\em presheaf pra monads}.

2. As it was pointed out in \cite{Leinster1} and \cite{Weber}  this framework fits well the free category monad over graphs and the free
$\o$-category monad over simple $\o$-graphs. I will elaborate on the first case as both cases are in a sense quite similar, well known and the first is simpler. In this case $p\Theta$ is the full subcategory of the category of graphs containing two graphs $[0]$ and $[1]$. The category of $T$-cardinals $\Theta_0$ is $\Delta_0$ and the category of $T$-cardinal algebras $\Theta_T$ is $\Delta$. The free category monad $T$ on  $Set^{\Theta^{op}_0}$ is pra
and the left adjoint $L_T$ to the functor $T_1:Set^{\Theta^{op}_0} \lra Set^{\Theta^{op}_0}\da T(1)$  can be described explicitly\footnote{In Example 2.5 of \cite{Weber} the functor $L_T$ is correctly described in words but $L_T$ is not given by the left Kan extension contrary to what was claimed there. I will came back to this point later.}. We shall sketch this definition to show the role of $\Delta_0$ in it. In many-to-one case the role of  $\Delta_0$ will be taken by $\ofs$.

Let $(G,|-|)$ be an object of the slice category of graphs $Graph\da T(1)$.  Thus we have a pair of function $d,c: E\ra V$ from the set edges to the set of vertices and a function $|-|: E\ra N$ from the set of edges  to the set of natural numbers. We define the diagram
\[ \Gamma_{(G,|-|)} : \tilde{G} \lra \Delta_0 \lra Graph \]
whose colimit is $L_T(G,|-|)$. The second functor is the usual embedding.
The set $\tilde{V}$ of vertices of $\tilde{G}$ contains both vertices and edges of $G$ as disjoint sets.  The set $\tilde{E}$ of edges of $\tilde{G}$ has two edges \[ d(e)\stackrel{s_e}{\lra}e,\hskip 5mm c(e)\stackrel{t_e}{\lra}e\]
 for each edge $e\in E$ with the domain and codomain as displayed. The functor $\tilde{G}\lra \Delta_0$ sends vertices from $V$ to $[0]$
 and the vertex $e\in E\subset \tilde{V}$ to the linear graph $[|e|]$. Moreover it sends the edges $s_e$ and $t_e$ to
\[ \bd^0_{[|e|]}:{[0]\lra[|e|]}, \mbox{ and}\;\;\; \bc^0_{[|e|]}: {[0]\lra[|e|]} \]
respectively. In particular, here and in all the other cases considered to get the formula for $L_T$ we use the domain and the codomain maps.

3. The case of positive-to-one computads also fits this setup and seems to be new.
The category of positive face structures $\fs$ is both the category of shapes for many-to-one computads $\Shape^{+/1}$
and the category of $T^{+/1}$-cardinals for the free $\o$-category monad $T^{+/1}$ on positive-to-one computads $\pComp$. Its image in $\oC$
is the category of $T^{+/1}$-cardinal algebras. The left adjoint $L_{T^{+/1}}$ to the $T^{+/1}_1$ functor can be described
much as in the previous case.

4. For the many-to-one computads the above setup does not seem to be sufficient. This is the first case where
the category of many-to-one shapes  $\mShape$ exists but it is not a full subcategory of the category $\mComp$ of many-to-one computads.
The category  $\mShape$ is equivalent to the category $\ofs$ of ordered face structures and monotone maps. But the category
of $T^{m/1}$-cardinals,  for the free $\o$-category monad $T^{m/1}$ on many-to-one computads $\mComp$  is equivalent to the category $\ofsl$ of ordered face structures and local maps. The category $\ofso$ of $T^{m/1}$-cardinal algebras will be described in the next section. Note that in the previous examples
the categories of cardinals were GT-theories but this time only the category of shapes $\ofs$ is a GT-theory and the category $\ofsl$ is not.

5. Finally let me point out one non-example namely the category of all computads $\Comp$. It is still true, by a beautiful argument of V. Harnik \cite{H},
that $\oC$ is monadic over $\Comp$ via right adjoint to the inclusion functor. However $\Comp$ is not a presheaf category, c.f. \cite{MZ2}, the free
$\o$-category monad on $\Comp$ is not pra as can be easily shown using Proposition 2.6 from \cite{Weber}. This adds to the long list of reasons why we don't
get a good theory of weak categories when considering all possible shapes of cells.

\section{GT-theories and nerves}\label{setup}

I formulate the general setup to show where it modifies the previous one. After stating the abstract pattern
I shall make a case study on the many-to-one computads to show the usefulness of this approach.  However
I do not provide any general results concerning this abstract setup as I prefer to collect more than one true example
($\ofs$) before developing this theory any farther. The present approach is in a sense much more modest than the
one from previous section. We deal here exclusively with the  monads whose categories
of algebras are equivalent to the category of  strict $\o$-categories $\oC$ only (or its truncations). Moreover we want these monad
to be defined on various reflective in $\oC$ subcategories of the category of computads $\Comp$. On the other
hand taking advantage of this more specific situation one may hope to get a more convenient description of concrete cases as in case
of the category of many-to-one computads $\mComp$. Still a word why all sorts of nerves of strict $\o$-categories might be of interest.
In the presheaf approach to weak categories (as opposed to the algebraic approach) the weak categories are presheaves with some properties.
If we believe that strict $\o$-categories should be special cases of weak ones, we need to study various nerves of strict $\o$-categories as they will provide abundance of examples.

The setup consists of GT-theory $\Phi$ with a full subcategory $p\Phi$ so that the obvious
induced functor  \[ Mod_\otimes(\Phi^{op},Set)\lra Set^{p\Phi^{op}} \]
is an equivalence of categories. With this data we want to get the following.
\begin{enumerate}
  \item The generic model $\cG_\Phi:\Phi\ra\Phi_l$ induces an equivalence of categories
  \[ Mod_\otimes(\Phi^{op},Set)\simeq pLim(\Phi_l^{op},Set) \]
  where $pLim(\Phi_l^{op},Set)$ is the category of functors that preserves the principal limits, i.e. the canonical limits defined over
  the diagrams consisting objects from $p\Phi$ only.
  \item There is a dense embedding
  \[ \Phi\lra Mod_\otimes(\Phi^{op},Set) \]
  \[ S\mapsto \bar{S}=\Phi_l(-,S) \]
  \item  $\Phi$ induces a pra monad $(T,\eta,\mu)$ on $Mod_\otimes(\Phi^{op},Set)$.
  \item  There is an explicit formula for the left adjoint $L_T$ to $T_1$.
  \item Let $\Phi_T$ denote the full image of $\Phi$ in $Alg(T)$ and $i:\Phi\ra\Phi_T$ the embedding.
  Then pra monadic functor  $Alg(T)\lra Mod_\otimes(\Phi^{op},Set)$ is a pseudo-pullback of the presheaf
  pra monadic functor $i^*:Set^{{\Phi_T}^{op}}\lra  Set^{{\Phi}^{op}}$, i.e. we have a pseudo-pullback
   \begin{center}
\xext=1000 \yext=500
\begin{picture}(\xext,\yext)(\xoff,\yoff)
 \setsqparms[1`1`1`1;1000`500]
 \putsquare(0,00)[Alg(T)`Set^{\Phi^{op}_T}` Mod_\otimes(\Phi^{op},Set)`Set^{\Phi^{op}_l}; `U`i^*`]
 \put(700,150){\makebox(100,100){$\cong$}}
 \end{picture}
\end{center}
  \item In particular, the  image of the full nerve functor $Alg(T)$ in  $Set^{{\Phi_T}^{op}}$
  consists of those functors that send $\otimes$-squares in $\Phi_T$ to pullbacks.
  \item The category of $T$-algebras $Alg(T)$ is equivalent to $\oC$.
\end{enumerate}

Now I will show how all this can be produced from the GT-theory $\ofs$ playing the role of $\Phi$,
together with its full subcategory of principal  ordered face structures $\pfs$ playing the role of $p\Phi$.

 The generic model of $\ofs$ is the inclusion functor $\ofs^{op}\ra \ofsl^{op}$. Thus the first two points are clear.
 In order to define the monad $T$ on $Mod_\otimes(\Phi^{op},Set)$ it is convenient to have already the full image of the
 category $\ofs$ in  $Alg(T)$ defined. This category $\ofso$, playing the role of $\Phi_T$, can be defined directly from
 $\ofs$ alone.  The objects of $\ofs_\o$ are the objects of  $\ofs$. A morphism
$\xi: R\ra S$ is  $\ofs_\o$ is an {\em $\o$-map} that is a transformation between presheaves $\xi : \ofsl(-,R)\lra \ofsl(-,S)$
which associate to a morphism $a: V\ra R$ in  $\ofsl(-,R)$ a morphism $\xi_a: V_a\ra S$ in  $\ofsl(-,S)$ so that
\begin{enumerate}
  \item $dim(V_a) \leq dim(V)$;
  \item $\xi(a\circ \bd^{(k)}_V)= \xi(a)\circ \bd^{(k)}_{V_a})$,\hskip 5mm $\xi(a\circ \bc^{(k)}_V)= \xi(a)\circ \bc^{(k)}_{V_a})$,
  \item if $V=V^1\otimes_kV^2$ then $V_a=V^1_{a\circ \kappa^1_{V^1}}\otimes_k V^2_{a\circ \kappa^2_{V^2}}$ and moreover
  \[ \xi(a\circ \kappa^1)= \xi(a) \circ \bar{\kappa}^1, \hskip 15mm
   \xi(a\circ \kappa^2)= \xi(a) \circ \bar{\kappa}^2 \]
   where $\kappa^i=\kappa^i_{V^i}$,  $\bar{\kappa}^i=\kappa^i_{V^i_{a\circ \kappa^i}}$ for $i=1,2$.
\end{enumerate}
We identify two $\o$-maps $\xi$ and $\xi'$ iff for every $a:V\ra R$ there is an isomorphism $\sigma_a $ making the triangle
\begin{center}
\xext=300 \yext=300
\begin{picture}(\xext,\yext)(\xoff,\yoff)
 \settriparms[1`1`1;300]
 \putVtriangle(0,0)[V_a`V'_a`S;\sigma_a`\xi_a`\xi'_a]
\end{picture}
\end{center}
commute.

Let  $\xi: R\ra S$ be an $\o$-map.  $\xi$ is an {\em inner $\o$-map} iff $\xi(1_R)=1_S$.  $\xi$ is a {\em monotone $\o$-map}
iff $\xi(1_R)$  a monotone morphism.
Its is easy to see that $\o$-maps, inner $\o$-maps and monotone $\o$-maps do compose. We have categories $\ofso$ of ordered face structures
and $\o$-map and $\ofsm$ of ordered face structures
and monotone $\o$-map. We write $f:R\ira S$ to indicate that $f$ is an inner map.
Clearly the $\o$-maps and the monotone $\o$-maps do compose. Thus we have categories $\ofso$ of  ordered face structures
and $\o$-maps, and $\ofsm$ of  ordered face structures
and monotone $\o$-maps. We have an embedding
\[ \iota_\o:\ofsl \lra \ofso \]
sending the local map $f: R\ra S$ to the $\o$-map $\iota_\o(f): R\ra S$ such that for $h:V\ra R$ in  $\ofsl(-,R)$ we have $\iota_\o(f)(h)=f\circ h$.
Clearly it restricts to the embedding \[ \iota_\mu:\ofs \lra \ofsm \]
Thus we have a commuting square of categories and functors
 \begin{center}
\xext=700 \yext=550
\begin{picture}(\xext,\yext)(\xoff,\yoff)
 \setsqparms[1`1`1`1;700`400]
 \putsquare(0,50)[\ofs`\ofsl`\ofsm`\ofso;\cG_\ofs`\iota_\mu`\iota_\o`\cG_{\ofsm}]
 \end{picture}
\end{center}
Note that the generic  models $\cG_\ofs$ and $\cG_{\ofsm}$ can be  considered as a process symmetrization of the tensor as
the pushouts tend to be more symmetric than tensors. Thus one cannot expect these functors to be full even on isomorphisms.
On the other hand, both $\iota_\mu$, $\iota_\o$ are full on isomorphisms and they fall under the scheme of the Example 1 from
the previous section. In particular, the composition functors induced by them $Set^{\ofsm^{op}}\lra Set^{\ofs^{op}}$ and
$Set^{\ofso^{op}}\lra Set^{\ofsl^{op}}$ are pra monadic.

Now the monad $T$ can be defined as follows. Let $X$ be a model in $Mod_\otimes(\Phi^{op},Set)$ and $k:R\ra S$ a monotone morphism.
For $S\in \ofs$ we put
\[ T(X)(S)= \coprod_{g:S\ira S'} X(S') \]
where the coproduct is taken over all (up to iso) inner $\o$-maps $g$ with the domain $S$. For a monotone morphism $k:R\ra S$,
the function $T(X)(k)$ is so defined that the square
 \begin{center}
\xext=1400 \yext=850
\begin{picture}(\xext,\yext)(\xoff,\yoff)
 \setsqparms[1`1`1`1;1700`600]
 \putsquare(0,100)[\;\;\;\;\;\;\;\;\;\;\;\;\;\;\;\;\;\;\;\;{T(X)(R)= \coprod_{f:R\sira R'} X(R')}`
 X(R')`\;\;\;\;\;\;\;\;\;\;\;\;\;\;\;\;\;\;\;\;{T(X)(S)= \coprod_{g:S\sira S'} X(S')}` X(S');
 \kappa^{X,R}_{f}`T(X)(k)`X(\bar{k})`\kappa^{X,S}_{g}]
 \end{picture}
\end{center}
commutes, for every inner $\o$-map $g: S\ira S'$ so that $f$ an inner $\o$-map and $\bar{k}$ is a monotone morphism and the square
 \begin{center}
\xext=400 \yext=450
\begin{picture}(\xext,\yext)(\xoff,\yoff)
 \setsqparms[1`1`1`1;400`300]
 \putsquare(0,50)[R`S`R`R';f`k`\bar{k}`g]
 \end{picture}
\end{center}
commutes. Clearly both $f$ and $\bar{k}$ are unique up an isomorphism. For the natural transformation $\alpha: X \ra Y$ in
$Mod_\otimes(\Phi^{op},Set)$ we define the natural transformation $T(\alpha)$ so that the square
 \begin{center}
\xext=1000 \yext=620
\begin{picture}(\xext,\yext)(\xoff,\yoff)
 \setsqparms[1`-1`-1`1;1000`500]
 \putsquare(0,20)[T(X)(S)`T(Y)(S)`X(S')`Y(S');T(\alpha)_S`\kappa^{X,S}_{g}`\kappa^{Y,S}_{g}`\alpha_{S'}]
 \end{picture}
\end{center}
commutes, for any $S\in \ofs$ and any inner $\o$-map $g:S\ra S'$. This ends the definition of $T$.

For $X\in Mod_\otimes(\Phi^{op},Set)$ the unit
\[ \eta_X: X\ra T(X) \]
is so define that for any $S\in \ofs$, we have $(\eta_X)_S=\kappa^{X,S}_{1_S}$.
The natural transformation $\eta$ is cartesian if for any morphism  $k:R\ra S$ in $\ofs$ the square
 \begin{center}
\xext=1200 \yext=650
\begin{picture}(\xext,\yext)(\xoff,\yoff)
 \setsqparms[1`-1`-1`1;1200`500]
 \putsquare(0,50)[X(R)`\coprod_{f:R\sira R'}X(R')`X(S)`\coprod_{f:S\sira S'}X(S');\kappa^{X,R}_{1_R}`X(k)`T(X)(k)`\kappa^{X,S}_{1_S}]
 \end{picture}
\end{center}
is a pullback. And this is so since the only way to complete the pair of morphisms
 \begin{center}
\xext=400 \yext=400
\begin{picture}(\xext,\yext)(\xoff,\yoff)
 \setsqparms[1`1`0`0;400`300]
 \putsquare(0,0)[R`R`S`;1_R`k``]
 \end{picture}
\end{center}
with $k$ being monotone morphism to a commuting square in $\ofso$
 \begin{center}
\xext=400 \yext=400
\begin{picture}(\xext,\yext)(\xoff,\yoff)
 \setsqparms[1`1`1`1;400`300]
 \putsquare(0,0)[R`R`S`S;1_R`k`k'`f]
 \end{picture}
\end{center}
with $f$ an inner $\o$-map and $k'$ a monotone morphism is by taking $k'=k$ and $f=1_S$.

For $X\in Mod_\otimes(\Phi^{op},Set)$ the multiplication
\[ \mu_X: T^2(X)\ra T(X) \]
is so define that for any $R\in \ofs$ and a pair of inner $\o$-maps $f:R\ra R'$ and $g:R'\ra R''$ the triangle
\begin{center}
\xext=1500 \yext=750
\begin{picture}(\xext,\yext)(\xoff,\yoff)
\putmorphism(0,700)(1,0)[T^2(X)(R)`T(X)(R)`(\mu_X)_S]{1500}{1}a
\putmorphism(0,580)(1,0)[ \parallel`\parallel`]{1500}{0}a
\putmorphism(0,480)(1,0)[\coprod_{f:R\sira R'} \coprod_{g:R'\sira R''}X(R'')`\coprod_{h:R\sira R''} X(R'')`]{1500}{0}a

\put(600,120){\vector(-2,1){500}}
\put(900,120){\vector(2,1){500}}

\put(700,0){\makebox(100,100){$X(R'')$}}

\put(130,150){\makebox(100,100){$\kappa^{X,S}_{f,g}$}}

\put(1200,150){\makebox(100,100){$\kappa^{X,S}_{g\circ f}$}}
\end{picture}
\end{center}
commutes. One can check that the natural transformation $\mu$ is cartesian if the commuting diagram
 \begin{center}
\xext=900 \yext=780
\begin{picture}(\xext,\yext)(\xoff,\yoff)
 \setsqparms[1`1`0`0;450`300]
 \putsquare(0,250)[R`R'`S`;f`k``]
 \setsqparms[1`0`1`0;450`300]
 \putsquare(450,250)[\phantom{R'}`R''`\phantom{S'}`S'';g``\check{k}`]

 \put(750,80){\vector(1,1){100}}
 \put(100,80){\line(-1,1){100}}
\put(100,80){\line(1,0){650}}

\put(400,0){\makebox(50,50){$h'$}}
 \put(750,710){\vector(1,-1){100}}
 \put(100,710){\line(-1,-1){100}}
\put(100,710){\line(1,0){650}}

\put(400,730){\makebox(50,50){$h$}}
 \end{picture}
\end{center}
in $\ofso$, with $f,g,h, h'$ inner $\o$-maps and $k,\check{k }$ monotone morphisms can be completed in a unique (up to an iso) way
to the commuting diagram
 \begin{center}
\xext=900 \yext=780
\begin{picture}(\xext,\yext)(\xoff,\yoff)
 \setsqparms[1`1`1`1;450`300]
 \putsquare(0,250)[R`R'`S`S';f`k``f']
 \setsqparms[1`0`1`1;450`300]
 \putsquare(450,250)[\phantom{R'}`R''`\phantom{S'}`S'';g`\bar{k}`\check{k}`g']

 \put(750,80){\vector(1,1){100}}
 \put(100,80){\line(-1,1){100}}
\put(100,80){\line(1,0){650}}

\put(400,0){\makebox(50,50){$h'$}}
 \put(750,710){\vector(1,-1){100}}
 \put(100,710){\line(-1,-1){100}}
\put(100,710){\line(1,0){650}}

\put(400,730){\makebox(50,50){$h$}}
 \end{picture}
\end{center}
with $f',g'$ inner $\o$-maps and $\bar{k}$ a monotone morphism. This is true as we can take as $S'$ the pushout of $f$ along $k$ in $\ofso$.
This property determines $S'$ up to isomorphism in $\ofso$ only. Then if we require additionally the lower order $<^\sim$ in $S'$ to be so
defined that $\bar{k}$ is monotone and $f'$ is inner, i.e. that all the left hand square is in $\ofsm$, then such an ordered face structure
$S'$ exists and is determined uniquely. Thus $\mu$ is cartesian indeed.

To show that  $(T,\eta,\mu)$ is a pra monad, it remains to show that the functor $T$ is a parametric right adjoint functor.
First let us describe $T(1)$. We can assume that for $R\in \ofs$
\[ T(1)(R)= \{ f:R\ira R' : f \in \ofs,\; f\, \mbox{inner} \}\cong \coprod_{ f:R\ira R'} 1(R') \]
Moreover, for $k:R\ra S$ in $\ofs$ and $g:S\ra S'$ an inner $\o$-map we have $T(1)(k)(g)=f$ if we have a square in $\ofso$
 \begin{center}
\xext=400 \yext=400
\begin{picture}(\xext,\yext)(\xoff,\yoff)
 \setsqparms[1`1`1`1;400`300]
 \putsquare(0,50)[R`R'`S`S';f`k`k'`g]
 \end{picture}
\end{center}
with $f$ an inner $\o$-map and $k'$ monotone morphism.
\[ T_1 :  Mod_\otimes(\Phi^{op},Set) \lra Mod_\otimes(\Phi^{op},Set) \da T(1) \]
\[ X\;\;\;\;\;\mapsto \;\;\;\;\; T(!_X):T(X)\lra T(1) \]
 has a left adjoint
 \[ L_T :  Mod_\otimes(\Phi^{op},Set)\da T(1) \lra Mod_\otimes(\Phi^{op},Set) \]
The construction of $L_T$ we give below is very similar in spirit to the construction from \mbox{Example 1} in Section \ref{pra monads}.
Let $|-|: Z\lra T(1)$ be an object of $Mod_\otimes(\Phi^{op},Set) \da T(1)$. We shall define a diagram
\[ \Gamma_{(Z,|-|)} : \cD_{(Z,|-|)} \lra Mod_\otimes(\Phi^{op},Set) \]
whose colimit is $L_T(Z,|-|)$. The set of objects of $\cD_{(Z,|-|)}$ is $\coprod_{S\in \ofs} Z(S)$, where as usual we take the coproduct over
isomorphism classes of objects of $\ofs$.   In other words an object of  $\cD_{(Z,|-|)}$ is a  pair $(y,S)$ so that $S\in \ofs$ $y\in Z(S)$.
We identify to such pairs $(y,S)$ and $(y',S')$ if there is a monotone isomorphism $f:S\ra S'$ such that $Z(h)(y')=y$. For any object $(y,S)$
and $k< dim(S)$ there are two arrows $d^{(k)}_{y,S}$  and $c^{(k)}_{y,S}$ in $\cD_{(Z,|-|)}$ with codomain $(y,S)$.  To describe the domains of
$d^{(k)}_{y,S}$  and $c^{(k)}_{y,S}$ let us note that since $|y|:S\ra R$ is an inner $\o$-map we can form a diagram
 \begin{center}
\xext=1600 \yext=550
\begin{picture}(\xext,\yext)(\xoff,\yoff)
 \setsqparms[1`1`1`1;800`400]
 \putsquare(0,50)[\bd^{(k)}S`S`\bd^{(k)}R`R;\bd^{(k)}_S`d^{(k)}(|y|)``\bd^{(k)}_R]
  \setsqparms[1`0`1`1;800`400]
 \putsquare(800,50)[\phantom{S}`\bc^{(k)}S`\phantom{R}`\bd^{(k)}R;\bc^{(k)}_S`|y|`c^{(k)}(|y|)`\bc^{(k)}_R]
 \end{picture}
\end{center}
with $\bd^{(k)}_S$, $\bc^{(k)}_S$,  $\bd^{(k)}_R$, $\bc^{(k)}_R$, monotone and $d^{(k)}(|y|)$, $c^{(k)}(|y|)$ inner $\o$-maps.
Then $d^{(k)}_{y,S}$  and $c^{(k)}_{y,S}$ have the domains as displayed in the diagram:
\begin{center}
\xext=2000 \yext=180
\begin{picture}(\xext,\yext)(\xoff,\yoff)
\putmorphism(0,20)(1,0)[(Z(\bd^{(k)}_S)(y),\bd^{(k)}S)`(y,S)`d^{(k)}_{(y,S)}]{1000}{1}a
\putmorphism(1000,20)(1,0)[\phantom{(y,S)}`(Z(\bc^{(k)}_S)(y),\bc^{(k)}S)`c^{(k)}_{(y,S)}]{1000}{-1}a
\end{picture}
\end{center}
We have a projection functor
\begin{center}
\xext=1800 \yext=200
\begin{picture}(\xext,\yext)(\xoff,\yoff)
\putmorphism(100,20)(1,0)[\cD_{(Z,|-|)}`\ofs`\pi_{(Z,|-|)}]{1700}{1}a
\end{picture}
\end{center}
which, in the notation as above with $|y|:S\ra R$, is given by
\begin{center}
\xext=3000 \yext=650
\begin{picture}(\xext,\yext)(\xoff,\yoff)
 \settriparms[-1`-1`0;500]
 \putAtriangle(300,50)[(y,S)`(Z(\bd^{(k)}_S)(y),\bd^{(k)}S)`(Z(\bc^{(k)}_S)(y),\bc^{(k)}S);d^{(k)}_{(y,S)}`c^{(k)}_{(y,S)}`]
 \settriparms[-1`-1`0;400]
 \putAtriangle(2100,150)[R`\bd^{(k)}R`\bc^{(k)}R;\bd^{(k)}_R`\bc^{(k)}_R`]

  \put(1500,300){\vector(1,0){300}}
 \put(1500,260){\line(0,1){80}}
\end{picture}
\end{center}
so that the composition of $\pi_{(Z,|-|)}$ with the embedding $\ofs \lra Mod_\otimes(\Phi^{op},Set)$ is the required diagram $\Gamma_{(Z,|-|)}$
for which we have $L_T(Z,|-|) = Colim \Gamma_{(Z,|-|)}$.

Next I will describe explicitly the adjunction:
\begin{center}
\xext=1200 \yext=220
\begin{picture}(\xext,\yext)(\xoff,\yoff)
\putmorphism(0,100)(1,0)[Mod_\otimes(\Phi^{op},Set)`\oC`]{1200}{0}a
\putmorphism(0,160)(1,0)[\phantom{Mod_\otimes(\Phi^{op},Set)}`\phantom{\oC}`F^{m/1}]{1200}{1}a
\putmorphism(0,40)(1,0)[\phantom{Mod_\otimes(\Phi^{op},Set)}`\phantom{\oC}`U^{m/1}]{1200}{-1}b
\end{picture}
\end{center}
For $X\in Mod_\otimes(\Phi^{op},Set)$ the set of $n$-cells is give by the coproduct
\[ F^{m/1}(X)_n = \coprod_{S\in \ofs_n} \coprod_{S\sira R} X(R) \]
The $k$-domain map of $n$-cells
operation is the unique morphism that for any $S\in\ofs_n$, and $f:S\ra R$ inner $\o$-map makes the square
  \begin{center}
\xext=1600 \yext=650
\begin{picture}(\xext,\yext)(\xoff,\yoff)
 \setsqparms[1`1`1`1;1200`500]
 \putsquare(0,50)[F^{m/1}(X)_n`F^{m/1}(X)_k`X(S)`X(\bd^{(k)}S);d^{(k)}`\kappa^{X,S}_{n,f}`\kappa^{X,\bd^{(k)}S}_{k,\bd^{(k)}(f)}`X(\bd^{(k)}_S)]
 \end{picture}
\end{center}
commute. The $k$-codomain map of $n$-cells is defined analogously and the identity operation is the obvious embedding.
Then the set of $k$-composable $n$-cells is
\[  F^{m/1}(X)_{n,k,n} = \coprod_{S,S'\in \ofs_n,\; \bc^{(k)}S=\bd^{(k)}S'}\;\; \coprod_{f:S\sira R,\; f':S'\sira R'} X(R\otimes_k R') \]
The composition morphism is the unique morphism making the triangles
\begin{center}
\xext=1400 \yext=700
\begin{picture}(\xext,\yext)(\xoff,\yoff)
 \settriparms[1`-1`-1;600]
 \putVtriangle(0,0)[ F^{m/1}(X)_{n,k,n}` F^{m/1}(X)_n`X(R\otimes_k R');m_{n,k,n}`\kappa^{X,S,S'}_{n,k,f,f'}`\kappa^{X,S\otimes_kS'}_{n,f\otimes_kf'}]
\end{picture}
\end{center}
commute, for any inner $\o$-maps  $f:S\ra R$,  $f':S'\ra R'$ so that  $\bc^{(k)}S=\bd^{(k)}S'$.
The $\o$-map $f\otimes_kf':S\otimes_kS'\lra R\otimes_kR'$ is well defined as as both $f$ and $f'$ are inner $\o$-maps.
This ends the definition of the $\o$-category $F^{m/1}(X)$. The definition of $F^{m/1}(X)$ on morphism is obvious.

The right adjoint to $F^{m/1}$ is induced by the embedding $\varepsilon:\ofs\ra \oC$. For $C\in\oC$ and $S\in \ofs$ we have
\[ U^{m/1}(C)(S)= \oC (\varepsilon(S),C) \]
The adjunction $F^{m/1}\dashv  U^{m/1}$ induces the above pra monad $(T,\eta,\mu)$. Moreover, by the Harnik argument c.f. \cite{H}\footnote{This argument was presented in \cite{MZ1} for the case of positive-to-one computads but it works without changes for the many-to-one computads as well.}, $U^{m/1}$ is monadic. Then from the theory developed by T. Leinster and M. Weber the left square is a pseudo-pullback
\begin{center}
\xext=2000 \yext=850
\begin{picture}(\xext,\yext)(\xoff,\yoff)
 \setsqparms[1`1`1`1;1000`500]

\put(1700,750){\vector(1,-1){150}}
\put(300,750){\line(1,0){1400}}
\put(150,600){\line(1,1){150}}
\put(1000,770){\makebox(100,100){$\cN_\mu$}}

 \putsquare(0,50)[\oC`Set^{\ofso^{op}}`Mod_\otimes(\ofs^{op},Set)`Set^{\ofsl^{op}};\cN_\o`U^{m/1}`\iota_\o^*`]
 \setsqparms[1`1`1`1;1000`500]
 \putsquare(1000,50)[\phantom{Set^{\ofso^{op}}}`Set^{\ofsm^{op}}`\phantom{Set^{\ofsl^{op}}}`Set^{\ofs^{op}};``\iota_\mu^*`]
 \end{picture}
\end{center}
Thus the pra monadic functor $U^{m/1}$ is a pseudo-pullback of pra monadic functor $\iota_\o^*$.
As the left bottom functor, induced by the generic model $\ofs\ra \ofsl$,  is full and faithful so is the many-to-one nerve functor $\cN_\o$ whose essential  image contains those presheaves whose restriction to $\ofsl^{op}$ preserves principal limits. It is not true that the right hand square in the above diagram is a  pseudo-pullback but the outer square still is and  the composition of the bottom functors is full and faithful. Hence the pra monadic functor $U^{m/1}$ is a pseudo-pullback of pra monadic functor $\iota_\mu^*$ as well. Thus we have another full nerve functor $\cN_\mu$ whose essential image consists of those functors whose restriction to $\ofs^{op}$ are models of $\ofs$.


\begin{thebibliography}{HMP}

\bibitem[BD]{BaezDolan}
 \frenchspacing
J. Baez, J. Dolan,  {\em Higher-dimensional algebra II:
n-Categories and the algebra of opetopes}. Advances in Math. 135
(1998), pp. 145-206.

\bibitem[B]{Batanin} M. Batanin, {\em Monoidal globular categories as
a natural environment for the theory of weak n-categories}.
Advances in Math. 136 (1998), 39-103.

\bibitem[B]{Berger} M. Berger, {\em A cellular nerve for higher categories},
Advances in Mathematics, 169, (2002), pp. 118-175.

\bibitem[C]{Cheng} E.Cheng, {\em The category of opetopes and the category of opetopic sets}, Theory and Applications of Categories, Vol. 11, No, 16, (2003), pp 353-374.

\bibitem[KJBM]{KJBM} J. Kock, A. Joyal, M. Batanin, J-F. Mascari,
 Polynomial  functors and opetopes, (2007), (math.QA/0706.1033).

\bibitem[H]{H} \mbox{V. Harnik}, private communication.

\bibitem[HMZ]{HMZ} \mbox{V. Harnik, M. Makkai, M. Zawadowski,}
{\em Multitopic sets are the same as many-to-one computads}, preprint
(2002), pp. 1-183, avaiable at http://www.math.mcgill.ca/makkai/.

\bibitem[HMP]{HMP} C. Hermida, M. Makkai, J. Power, {\em On weak
higher dimensional categories I} Parts 1,2,3, J. Pure and Applied
Alg.  153 (2000), pp. 221-246, 157 (2001), pp. 247-277, 166
(2002), pp. 83-104.

\bibitem[J]{Joyal} A.\ Joyal, {\em Disks, Duality and $\Theta$-categories}. Preprint,
(1997).

\bibitem[Law]{Lawvere} F.W. Lawvere, {\em Functorial semantics of algebraic theories}, desertation.
Columbia Iniv., New York (1963), (summary in Proc. Nat. Acad. Sci 50, (1963), pp. 869-872).

\bibitem[Lei1]{Leinster0} T. Leinster, {\em Higher Operads, Higher Categories}. LOndon Math. Soc. Lecture Note Series. Cambridge University Press, Cambridge, 2004. (math.CT/0305049)

\bibitem[Lei1]{Leinster1} T. Leinster, {\em Nerves of Algebras}, talk at CT04 in Vancouver,
avaiable at http://www.maths.gla.ac.uk/~tl/vancouver/.

\bibitem[Lei2]{Leinster2} T. Leinster, {\em How I learned to Love the Nerves Functor},
An essey on The n-Category Caf\'{e}, Janurary, 2008.

\bibitem[Lin]{Linton} F.E.J. Linton, {\em Some aspects of Equational Categories},
Proc. of the Conference on Categorical Algebra, La Jolla 1965,
 Springer-Verlag, (1966), pp. 84-94.

\bibitem[M]{M} M. Makkai,
{The multitopic $\o$-category of all multitopic $\o$-categories}, preprint
(1999), pp. 1-67, avaiable at http://www.math.mcgill.ca/makkai/.

\bibitem[P]{Palm} T. Palm, {\em Dendrotopic sets}, Hopf algebras, and semiabelian categories, Fields Inst. Commun.
vol. 43 (2004), 411-461 AMS, Providence, RI.

\bibitem[MZ1]{MZ1}
M. Makkai and M. Zawadowski, {\em Disks and duality}. TAC 8(7),
2001, pp. 114-243.

\bibitem[MZ2]{MZ2} M. Makkai, M.Zawadowski,
{The category of 3-computads is not cartesian closed}, to Appear in J. of Pure and Applied Algebra.

\bibitem[W]{Weber}
M. Weber, {\em Familial 2-functors and parametric right adjoints}.
Theory and Applications of Categories, vol. 18,  (2007), pp. 665-732.

\bibitem[Z1]{Z1}
M. Zawadowski, {\em On positive face structures and
positive-to-one computads}. Preprint, (2006), pp. 1-77.

\bibitem[Z2]{Z2}
M. Zawadowski, {\em On ordered face structures and many-to-one
computads}. Preprint, (2007), pp. 1-95.
\end{thebibliography}
\end{document}